\documentclass[11pt]{amsart}
\usepackage{amssymb, pstricks,pstcol,pst-plot}

\definecolor{OrangeRed}{cmyk}{0,0.6,1,0}            
\definecolor{DarkBlue}{cmyk}{1,1,0,0.20}
\definecolor{Black}{cmyk}{0,0,0,1}
\definecolor{Violet}{cmyk}{0.79,0.88,0,0}
\definecolor{myblue}{rgb}{0.66,0.78,1.00}

\parskip=\smallskipamount
\numberwithin{equation}{section}

\newtheorem{theorem}{Theorem}[section]
\newtheorem{lemma}[theorem]{Lemma}
\newtheorem{corollary}[theorem]{Corollary}
\newtheorem{proposition}[theorem]{Proposition}

\theoremstyle{definition}
\newtheorem{definition}[theorem]{Definition}

\newtheorem{remark}[theorem]{Remark}

\newcommand{\B}{\mathbb{B}}
\newcommand{\C}{\mathbb{C}}

\newcommand{\N}{\mathbb{N}}

\newcommand{\Z}{\mathbb{Z}}
\renewcommand{\P}{\mathbb{P}}
\newcommand{\R}{\mathbb{R}}

\newcommand{\cA}{\mathcal{A}}

\newcommand{\cC}{\mathcal{C}}

\newcommand{\cL}{\mathcal{L}}

\newcommand{\cO}{\mathcal{O}}

\newcommand{\cZ}{\mathcal{Z}}
\newcommand{\cP}{\mathcal{P}}

\newcommand\wt{\widetilde}
\newcommand\spsh{strongly plurisubharmonic}
\newcommand\hra{\hookrightarrow}

\def\ss{\Subset}
\def\di{\partial}
\def\dibar{\bar\partial}
\def\bs{\backslash}

\def\e{\epsilon}

\def\l{\lambda}
\def\d{\delta}

\def\disc{\triangle}

\def\cn{{\mathbb C}^n}

\def\dim{{\rm dim}\,}    
    
\def\psh{plurisubharmonic}              
\def\nbd{neighborhood}

\begin{document}
\title{Holomorphic curves in complex spaces}
\author{Barbara Drinovec-Drnov\v sek \& Franc Forstneri\v c}
\address{Institute of Mathematics, Physics and Mechanics, 
University of Ljubljana, Jadranska 19, 1000 Ljubljana, Slovenia}
\email{barbara.drinovec@fmf.uni-lj.si}
\email{franc.forstneric@fmf.uni-lj.si}
\thanks{Research supported by grants P1-0291 and J1-6173, Republic of Slovenia.}

%
%

\date{March 9, 2006} 

\begin{abstract}
We study the existence of topologically closed complex cur\-ves  normalized 
by bordered Riemann surfaces in complex spaces. Our main result is that such curves 
abound in any noncompact complex space admitting an exhaustion function whose 
Levi form has at least two positive eigenvalues at every point outside 
a compact set, and this condition is essential. We also construct 
a Stein neighborhood basis of any compact complex curve with $\cC^2$ boundary 
in a  complex space.
\end{abstract}
\maketitle

{\small \rightline{\em  To Josip Globevnik for his 60th birthday}}

%
%
%
%
\section{Introduction}
Let $X$ be an irreducible (reduced, paracompact) complex space
of dimension $>1$.  For every topologically closed complex curve $C$
in $X$ we have a sequence of holomorphic maps
\[
	\{\C\P^1,\C,\disc\}\ni \wt D \to D \to C \hra X
\]
where $C\hra X$ is the inclusion, $D\to C$ is a normalization of $C$
by a Riemann surface $D$, and $\wt D\to D$ is a universal covering 
combined with a uniformization map. Here $\disc=\{z\in\C\colon |z|<1\}$.  Thus $C$ 
is the image of a generically one to one proper holomorphic map $D \to X$;
hence it is natural to ask which Riemann surfaces $D$ admit any proper 
holomorphic maps to a given complex space, and how plentiful are they. 
This question has been investigated most intensively for compact complex curves
which form a part of the {\em Douady space} and of the 
{\em cycle space} of $X$ (\cite{Barlet}, \cite{Campana}, \cite{Douady}).

In this paper we obtain essentially optimal existence and approximation
results when $D$ is a {\em finite bordered Riemann surface},
i.e., a one dimen\-sional complex manifold with compact closure
$\bar D=D\cup bD$ whose boundary $bD$ consists of finitely many closed 
Jordan curves; such $D$ is uniformized by the disc $\disc$.
The existence of a proper holomorphic map $D\to X$ implies  that 
$X$ is noncompact, but additional conditions are needed in general 
since there exist open  complex manifolds without 
any topologically closed complex curves; an example is obtained by 
removing a point from a compact complex manifold which admits 
no closed complex curves (a condition satisfied e.g.\ by certain  
complex tori of dimension $>1$).

We begin by a brief survey of the known results.
Every open Riemann surface admits a proper holomorphic immersion in $\C^2$ 
and a proper holomorphic embedding in $\C^3$ \cite{Bishop}, \cite{Nar1}. 
Some open Riemann surfaces also embed in $\C^2$,
but it is unknown whether all of them do; impressive results on 
this subject have been obtained recently by E.\ Forn\ae ss Wold  
(\cite{Wo0}, \cite{Wo1}, \cite{Wo2}) where the reader can find references 
to older works on the subject. 

Turning to more general target spaces, 
we note that the Kobayashi hyper\-bolicity of $X$ 
excludes curves uniformized by $\C$ but imposes less restrictions 
on those uniformized by  the disc $\triangle$ \cite{Ko1}, \cite{Ko2}. 
There are other, less tangible obstructions:  Dor \cite{Dor2} found 
a bounded domain with non smooth boundary in $\C^n$ without any proper 
holomorphic images of $\disc$; even in  smoothly bounded (non pseudoconvex)
domains in $\C^n$ the union of images of all proper analytic discs 
can omit a nonempty open subset \cite{FGl}.
On the positive side, every point in a Stein manifold $X$
of dimension $>1$ is contained in the image of a proper holomorphic 
map $\disc \to X$ (Globevnik \cite{Glo}; see also \cite{Dor1}, \cite{BD1}, 
\cite{BD2}, \cite{BD3}, \cite{FGl}, \cite{FG2}, \cite{FGS}).
The same holds for discs in any connected complex manifold $X$ 
which is $q$-complete for some $q < \dim X$ \cite{BD3}. 
The first case of interest, inaccessible with 
the existing techniques, are Stein spaces with singularities.

Recall that a smooth function 
$\rho\colon X\to\R$ on a complex space $X$ is said to be {\em $q$-convex} 
on an open subset $U\subset X$ (in the sense of Andreotti-Grauert \cite{AG},
\cite[def.\ 1.4, p.\ 263]{Grauert2})
if there is a covering of $U$ by open sets $V_j \subset U$, biholomorphic
to closed analytic subsets of open sets $\Omega_j\subset\C^{n_j}$,
such that for each $j$ the restriction $\rho|_{V_j}$ admits an extension
$\wt\rho_j\colon \Omega_j\to\R$ whose Levi form $i\di\dibar \,\wt\rho_j$ 
has at most $q-1$ negative or zero eigenvalues at each point of $\Omega_j$.
The space $X$ is {\em $q$-complete}, resp.\ {\em $q$-convex}, if it admits 
a smooth exhaustion function $\rho\colon X\to\R$ which is $q$-convex on $X$, 
resp.\ on $\{x\in X\colon \rho(x)>c\}$ for some $c\in \R$.
A 1-complete complex space is just a Stein space, and a 1-convex space
is a proper modification of a Stein space. 
We denote by $X_{reg}$ (resp.\ by $X_{sing}$) the set of regular 
(resp.\ singular) points of $X$.

We are now ready to state our first main result; it is proved in \S 6.

%
%
%
%
\begin{theorem}
\label{main}
Let $X$ be an irreducible complex space of $\dim X>1$, and let
$\rho\colon X\to\R$ be a smooth exhaustion function which is $(n-1)$-convex 
on $X_c=\{x\in X\colon\rho(x)> c\}$ for some $c\in\R$.
Given a bordered Riemann surface $D$ and a $\cC^2$ map  
$f\colon \bar D\to X$ which is holomorphic in $D$ and satisfies 
$f(D)\not\subset X_{sing}$ and $f(bD) \subset X_c$, there is a sequence 
of proper holomorphic maps $g_{\nu} \colon D\to X$ homotopic to $f|_D$
and converging to $f$ uniformly on compacts in $D$ as $\nu\to\infty$. 
Given an integer $k\in\N$ and finitely many points $\{z_j\}\subset D$, 
each $g_\nu$ can be chosen to have the same $k$-jet as $f$ at each of the points $z_j$.
\end{theorem}

We now show by examples that the conditions in the above theorem are essentially optimal. 
The assumption on $\rho$ means that its Levi form has 
at least two positive eigenvalues at every point of $X_{c} = \{\rho>c\}$.
One positive eigenvalue does not suffice in view of Dor's example of a domain  
in $\C^n$ without any proper analytic discs \cite{Dor2} and the fact that
every domain in $\C^n$ is $n$-complete (\cite{GW}, \cite{Ohs}).
Necessity of the hypothesis
$f(D)\not\subset X_{sing}$ is seen by Proposition 3 in \cite{FW} 
(based on an example of Kaliman and Zaidenberg \cite{KZ}): An analytic disc 
contained in $X_{sing}$ may be forced to remain there under analytic 
perturbations, and it need not be 
approximable by proper holomorphic maps $\disc\to X$. 
The only possible improvement would be a reduction of 
the boundary regularity assumption on the initial map. 
If $D$ is a planar domain bounded by finitely many 
Jordan curves and $X$ is a manifold,  it suffices to 
assume that $f$ is continuous on $\bar D$ by appealing 
to \cite[Theorem 1.1.4]{Chak} in order to approximate $f$ by a 
more regular map.     

If $f\colon\bar D\to X$ in Theorem \ref{main} is generically injective  then
so is any proper holomorphic map $g_\nu \colon D\to X$
approximating $f$ sufficiently closely; its image $g_\nu(D)$ 
is then a closed complex curve in $X$ normalized by $D$. 
Assuming that $f(\bar D)\subset X_{reg}$ one can choose each $g_{\nu}$ 
to be an immersion, and even an embedding when $n\ge 3$.
Each map $g_\nu$ will be a locally uniform limit in $D$ of a sequence of $\cC^2$ maps 
$f_j\colon \bar D\to X$ which are holomorphic in $D$ and satisfy
\begin{equation}
\label{to-infinity}
	\lim_{j\to\infty} \inf\{\rho\circ f_j(z) \colon z\in bD\}  \to +\infty;
\end{equation}
that is, their boundaries $f_j(bD)$ tend to infinity in $X$. 
Embedding $\bar D$ as a domain in an open Riemann surface $S$, we can 
choose each $f_j$  to be holomorphic in open set $U_j\subset S$ 
containing $\bar D$.

Theorem \ref{main} also gives new information on {\em algebraic curves} 
in $(n-1)$-convex {\em quasi projective algebraic spaces} $X=Y\bs Z$, where 
$Y,Z\subset\C\P^N$ are closed complex ($=$algebraic) 
subvarieties in a complex projective space.
We embed our bordered Riemann surface $D$ as a domain with smooth 
real analytic boundary in its {\em double} $\hat S$, a compact 
Riemann surface obtained by gluing two copies of $\bar D$ along their 
boundaries \cite[p.\ 581]{BS}, \cite[p.\ 217]{Spr}.
There is a meromorphic embedding $\hat S\hra\C\P^3$ with poles outside of 
$\bar D$;  the subset $S \subset \hat S$ which is mapped to the affine part  
$\C^3 \subset \C\P^3$ is a smooth affine algebraic curve, and $D$ 
is Runge in $S$. A holomorphic map $f\colon U\to X$ 
from an  open set $U\subset S$ to a quasi projective algebraic space 
$X$ is said to be {\em Nash algebraic} (Nash \cite{Nash}) if the graph 
\[
	G_f=\{(z,f(z)) \in S\times X\colon z\in U\}
\]
is contained in a one dimensional algebraic subvariety of $S\times X$.

%
%
%
%
\begin{corollary}
\label{algebraic}
Let $X$ be an irreducible quasi projective algebraic space of $\dim X>1$,
and let $D\Subset S$ be a smoothly bounded Runge domain in an affine algebraic
curve $S$. Assume that $\rho\colon X\to\R$ and $f\colon \bar D \to X$ 
satisfy the hypotheses of Theorem \ref{main}.
Then there is a sequence of Nash algebraic maps $f_j\colon U_j\to X$ in open sets
$U_j\supset \bar D$ satisfying (\ref{to-infinity}) such that the sequence 
$f_j|_D$ converges to a proper holomorphic map $g\colon D\to X$. 
\end{corollary}

Corollary \ref{algebraic} is obtained by approximating each of the 
holomorphic maps $f_j\colon U_j\to X$, obtained in the proof 
of Theorem \ref{main}, uniformly on $\bar D$ by a Nash algebraic map,
appealing to theorems of Demailly, Lempert and Shiffman \cite[Theorem 1.1]{DLS} 
and Lempert \cite{Lempert}. Their results  give Nash algebraic 
approximations of any holomorphic map from an open Runge domain in an 
affine algebraic variety to a quasi projective algebraic space.
Of course $g$ can be chosen to also satisfy the additional 
properties in Theorem \ref{main}. If $\Gamma_j\subset S \times X$ 
is an algebraic curve containing the graph of the Nash algebraic map 
$f_j\colon U_j\to X$ then its projection $C_j\subset X$ under the map $(z,x)\to x$ 
is an algebraic curve in $X$ containing $f_j(U_j)$;
as $j\to \infty$, the domains $f_j(D) \subset C_j$ 
converge to the closed transcendental curve $g(D) \subset X$
while their boundaries $f_j(bD)$ leave any compact subset of~$X$. 

Corollary \ref{algebraic} applies for example to $X=\C\P^n \bs A$ 
where $A$ is a closed complex  submanifold of dimension 
$d\in\{[\frac{n+1}{2}], \ldots,n-1\}$. Indeed, 
$\C\P^n \bs A$ is then $(2(n-d)-1)$-complete by a 
result of Peternell \cite{Peternell} (improving an earlier 
result of Barth \cite{Barth}), and hence is $(n-1)$-complete if $n\le 2d$.

%
%
Another interesting and relevant example is due to Schneider \cite{Schneider}
who proved that for a compact complex manifold $X$ 
and a complex submanifold $A\subset X$ of codimension $q$
whose normal bundle $N_{A|X}$ is (Griffiths) positive
the complement $X\bs A$ is $q$-convex. Thus Theorem \ref{main}
furnishes closed complex curves in $X\bs A$ whenever 
$q\le \dim X-1$, which is equivalent to $\dim A\ge 1$. 
For further examples see Grauert \cite{Grauert2} and
Col\c toiu \cite{Col1}.

The following consequence of Theorem \ref{main} was proved in 
\cite{BD3} in the special case when $X_{sing}=\emptyset$ and $D=\disc$.  

%
%
%
%
\begin{corollary}
\label{complete-spaces}
Let $X$ be an irreducible $(n-1)$-complete complex space of dimension $n>1$, and
let $D$ be a bordered Riemann surface. Given a $\cC^2$ map $f\colon \bar D\to X$ 
which is holomorphic in $D$ and satisfies $f(D)\not\subset X_{sing}$, a positive
integer $k\in\N$ and finitely many points $\{z_j\} \subset D$, there is a sequence 
of proper holomorphic maps $g_\nu \colon D\to X$ converging to $f|_D$ 
uniformly on compacts in $D$ such that each $g_\nu$ has the same $k$-jets as $f$ 
at each of the points $z_j$. This holds in particular if $X$ is a Stein space.
\end{corollary}

Let $X$ be a complex manifold. The {\em Kobayashi-Royden pseudo\-norm} 
of a tangent vector $v\in T_x X$ is given by
\[
	\kappa_X(v)=\inf \left\{ \lambda>0 \colon 
	\exists f\colon \disc\to X\ {\rm holomorphic},\ f(0)=x,\ f'(0)=\lambda^{-1}v \right\}.
\]
The same quantity is obtained by using only maps which are holomorphic 
in small neighborhoods of $\bar \disc$ in $\C$. Corollary 
\ref{complete-spaces} implies:

%
%
%
%
\begin{corollary}
\label{Kobayashi}
If $X$ is an $(n-1)$-complete complex manifold of dimension $n>1$
then its infinitesimal Kobayashi-Royden pseudometric $\kappa_X$ is computable
in terms of {\em proper} holomorphic discs $f\colon \disc\to X$.
\end{corollary}

On a quasi projective algebraic manifold $X$, the pseudometric $\kappa_X$ 
and its integrated form, the Kobayashi pseudodistance, are also
computable by algebraic curves \cite[Corollary 1.2]{DLS}.

It is natural to inquire which homotopy classes of maps $D\to X$ from a 
bordered Riemann surface admit a proper holomorphic representative. 
Hyperbolicity properties of $X$ may impose a major obstruction 
on the existence of a holomorphic map in a given nontrivial homotopy class
(\cite{Ko1}, \cite{Ko2}, \cite{Eisenman}). The following opposite property 
is important in the Oka-Grauert theory: 

A complex manifold $X$ is said to enjoy  the 
{\em $m$-dimensional convex approximation property} $({\rm CAP}_m)$ 
if every holomorphic map $U\to X$ from an open set $U \subset \C^m$
can be approximated uniformly on any compact  convex set $K\subset U$ 
by entire maps $\C^m\to X$ \cite{ANN}.

\begin{corollary}
\label{cor2}
Let $X$ be an $(n-1)$-complete complex manifold of dimension $n>1$. 
If $X$ satisfies ${\rm CAP}_{n+1}$ then for every continuous map 
$f\colon D\to X$ from a  bordered Riemann surface $D$ there 
exists a proper holomorphic map $g\colon D\to X$ homotopic to $f$. 
If $f$ is holomorphic on a neighborhood of a compact subset $K\subset D$
then $g$ can be chosen to approximate $f$ as close as desired on $K$.
This holds in particular if $X=\C\P^n \bs A$ where $n\ge 4$ 
and $A\subset \C\P^n$ is a closed complex submanifold of 
dimension $d\in\{[\frac{n+1}{2}], \ldots,n-2\}$. 
\end{corollary}

\begin{proof}
We may assume that $\bar D=\{z\in S\colon v(z)\le 0\}$ 
where $S$ is an open Riemann surface and $v\colon S \to \R$
is a smooth function with $dv\ne 0$ on $bD=\{v=0\}$.
Choose numbers $c_0 <0< c_1$ close to $0$ such that $v$ has no critical 
values on $[c_0,c_1]$. Let $D_j=\{z\in S\colon v(z)<c_j\}$
for $j=0,1$. We may assume $K\subset D_0$. 
There is a homotopy of smooth maps 
$\tau_t\colon D_1\to D_1$ $(t\in [0,1])$ such that $\tau_0$ 
is the identity on $D_1$, $\tau_1(D_1)=D_0$, and for all
$t\in[0,1]$ we have $\tau_t(D)\subset D$ and $\tau_t$
equals the identity map near $K$. Set $\wt f= f\circ \tau_1\colon D_1\to X$.
Note that $\wt f|_D$ is homotopic to $f$ via the homotopy $f\circ \tau_t|_D$
$(t\in[0,1])$. 

By the main result of \cite{ANN} the ${\rm CAP}_{n+1}$ property 
of $X$ implies the existence of a holomorphic map $f_1 \colon D_1 \to X$ 
homotopic to $\wt f\colon D_1 \to X$. Then $f_1|_D$ is homotopic 
to $\wt f|_D$ and hence to $f$. Theorem \ref{main}, applied to the map 
$f_1|_{\bar D} \colon\bar D\to X$, furnishes a proper holomorphic map 
$g\colon D\to X$ homotopic to $f_1|_D$, and hence to $f$. 
In addition, $f_1$ and $g$ can be chosen to approximate
$f$ uniformly on $K$.

The last statement follows from the already mentioned fact that 
$\C\P^n \bs A$ is $(n-1)$-complete if $A$ is as in the statement
of the corollary (see \cite{Peternell}), and it enjoys ${\rm CAP}_{m}$ 
for all $m\in\N$ provided that $\dim A\le n-2$ \cite{ANN}. 
\end{proof}

By \cite{ANN} and \cite{FFourier} the property ${\rm CAP}=\cap_{m=1}^\infty {\rm CAP}_m$
of a complex manifold $X$ is equivalent to the classical {\em Oka property}
concerning the existence and the homotopy classification of holomorphic maps
from Stein manifolds to $X$. Examples in \cite{Gromov} and \cite{ANN}
show that Corollary \ref{cor2} fails in general if $X$ does not enjoy CAP,
and the most one can expect is to find a proper map $D\to X$ 
in the given homotopy class which is holomorphic with respect to 
some complex structure on the smooth 2-surface $D$.
This indeed follows by combining Theorem \ref{main} with 
a very special case of the main result in \cite{FS1}.

\begin{corollary}
\label{cor3}
Let $X$ be a $(n-1)$-complete complex manifold of dimension $n>1$
and let $\bar D$ be a compact, connected, oriented real surface with boundary.
For every continuous map $f\colon D\to X$ there exist a 
complex structure $J$ on $D$ and a proper $J$-holomorphic map 
$g\colon D\to X$ which is homotopic to $f$. 
\end{corollary}

Another result of independent interest is Theorem \ref{Stein-nbds}
to the effect that a compact complex curve with $\cC^2$ boundary 
in a complex space admits a basis of open Stein neighborhoods.  
The following special case is proved in~\S 2.

%
%
%
%
\begin{theorem}
\label{Euclidean-nbds}
Let $X$ be an $n$-dimensional complex manifold. 
If $D$ is a relatively compact smoothly bounded domain in an open 
Riemann surface $S$ and $f\colon\bar D\hra X$ is a $\cC^2$ embedding which is 
holomorphic in $D$ then $f(\bar D)$ has a basis of open Stein neighborhoods in $X$
which are biholomorphic to open neighborhoods of $\bar D\times\{0\}^{n-1}$ 
in $S\times \C^{n-1}$.
In particular, if $D$ is a smoothly bounded planar domain 
then $f(\bar D)$ has a basis of open Stein neighborhoods in $X$
which are biholomorphic to domains in $\C^{n}$.
\end{theorem}

Royden showed in \cite{Royden} that for any holomorphically embedded polydisc 
$f\colon \disc^k\hra X$ in a complex manifold $X$ and for any $r<1$
the smaller polydisc $f(r\disc^k) \subset X$ admits open neighborhoods in $X$ 
biholomorphic to $\disc^n$ with $n=\dim X$. We have the analogous result 
for closed analytic discs, showing that they have no appreciation whatsoever 
of their surroundings.

\begin{corollary}
\label{discs}
Let $X$ be an $n$-dimensional complex manifold.
For every $\cC^2$ embedding $f\colon\bar\disc \hra X$ which is 
holomorphic in $\disc$ the image $f(\bar \disc)$ has a basis
of open neighborhoods in $X$ which are biholomorphic to $\disc^n$.
\end{corollary}

These and related result are used to obtain new holomorphic approximation theorems
(Corollary \ref{approximation} and Theorem \ref{approximation2}).

%
%
%
%
\smallskip
{\em Outline of proof of Theorem \ref{main}}.
Theorem \ref{main} is proved in \S 6 after developing the necessary tools  in \S 2--\S 5. 
We begin by perturbing the initial map $f\colon \bar D\to X$ to a new map 
for which $f(bD)\subset  X_{reg}$ (Theorem \ref{approximation2}).
The rest of the construction is done in such a way that the image 
of $bD$ remains in the regular part of $X$.
A proper holomorphic map $g\colon D\to X$ 
is obtained as a limit $g=\lim_{j\to\infty} f_j|_D$ of a sequence of
$\cC^2$ maps $f_j\colon \bar D\to X$ which are holomorphic in $D$ such that 
the boundaries $f_j(bD)$ converge to infinity. 

Our local method of lifting the boundary $f(bD)$  is similar 
to the one used (in the special case $D=\disc$) in earlier papers on the subject 
(\cite{Dor1}, \cite{BD1}, \cite{BD2}, \cite{FGl}, \cite{FG2}, \cite{Glo}).  
Since the Levi form $\cL_\rho$ is assumed to have at least two positive eigenvalues 
at every point of $f(bD)$, we get at least one positive eigenvalue in a direction 
tangential to the level set of $\rho$ at each point $f(z)$, $z\in bD$;
this gives a small analytic disc in $X$, tangential to the 
level set of $\rho$ at $f(z)$, along which $\rho$ increases quadratically. 
By solving a certain Riemann-Hilbert boundary value problem we obtain a 
local holomorphic map whose boundary values on the relevant part of $bD$ 
are close to the boundaries of these discs, and hence $\rho\circ f$ has 
increased there. (One positive eigenvalue of $\cL_\rho$ does not 
suffice since the corresponding eigenvector may be transverse to 
the level set of $\rho$ and cannot be used in the construction.)

To globalize the construction we develop a new method of patching holomorphic maps
by improving a technique from the recent work of the second 
author on localization of the Oka principle \cite{ANN}. 
We embed a given map $f\colon \bar D\to X$  into a {\em spray of maps}, i.e.,
a family of maps $f_t\colon \bar D\to X$ depending holomorphically 
on the parameter $t$ in a Euclidean space 
and satisfying a certain submersivity property 
(dominability) outside of an exceptional subvariety. 
The local modification method explained above gives a new spray 
near a part of the boundary $bD$; by insuring that the two sprays are sufficiently 
close to each other on the intersection of their domains $\overline{D_0 \cap D}_1$, 
we patch them into a new spray  over $\overline {D_0\cup D}_1$ 
(Proposition \ref{gluing-sprays}). This is accomplished 
by finding a fiberwise biholomorphic transition map between them and decomposing 
it into a pair of maps over $\bar D_0$ resp.\ $\bar D_1$ which are used 
to correct the two sprays so as to make them agree over $\overline{ D_0\cap D}_1$. 

The main step, namely a decomposition of the transition map 
(Theorem \ref{splitting}), is achieved by a rapidly convergent iteration.
This result generalizes the classical Cartan's lemma to non linear maps, 
with $\cC^r$ estimates up to the boundary. Unlike in  \cite[Lemma 2.1]{ANN}, 
the base domains don't shrink in our present construction --- this is not allowed since 
all action in the construction of proper maps takes place at the boundary.

Our method of gluing sprays is also useful in proving holomorphic
approximation theorems (see Theorem \ref{approximation2} below).

One of the difficult problems in earlier papers has been to avoid
running into a critical point of the given exhaustion function 
$\rho\colon X\to\R$. For Stein manifolds this problem was  
solved by Globevnik \cite{Glo}. Here we apply an alternative method 
from \cite{ACTA} and cross each critical level by using a 
different function constructed especially for this purpose. 

We believe that the methods developed in this paper will be applicable in 
other problems involving holomorphic maps.  With this in mind, 
many of the new technical tools are obtained in the more general context 
of strongly pseudoconvex domains in Stein manifolds.

%
%
%
%
%
%
\section{Stein neighborhoods of smoothly bounded complex curves}
Let $(X,\cO_X)$ be a complex space. 
We denote by $\cO(X)$ the algebra of all holomorphic functions 
on $X$, endowed with the compact-open topology.  A compact subset 
$K$ of $X$ is said to be {\em $\cO(X)$-convex} if for any point $p\in X\backslash K$ 
there exists $f\in \cO(X)$ with $|f(p)|>\sup_K |f|$. If $X$ is Stein
and $K$ is contained in a closed complex subvariety $X'$ of $X$ 
then $K$ is $\cO(X')$-convex if and only if it is $\cO(X)$-convex. 
(For Stein spaces we refer to \cite{GR} and \cite{Ho}.) 

We will say that a compact set $A$ in a complex space $X$ 
is a {\em complex curve with $\cC^r$ boundary $bA$ in $X$} if 
\begin{itemize}
\item[(i)] $A\bs bA$ is a closed, purely one dimensional 
complex subvariety of $X\bs bA$ without  compact irreducible 
components, and
\item[(ii)] every point $p\in bA$ has an open neighborhood $V\subset X$
and a biholomorphic map $\phi\colon V\to V'\subset \Omega\subset \C^N$
onto a closed complex subvariety $V'$ in an open subset $\Omega \Subset\C^N$
such that $\phi(A\cap V)$ is a one dimensional complex submanifold of $\Omega$ with 
$\cC^r$ boundary $\phi(bA \cap V)$. 
\end{itemize}

Note that $bA$ consists of finitely many closed Jordan curves and 
has no isolated points, but it may contain some singular points of $X$.

%
%
%
%
\begin{theorem}
\label{Stein-nbds}
Let $A$ be a compact complex curve with $\cC^2$ boundary in a complex space $X$.
Let $K$ be a compact $\cO(\Omega)$-convex set in a Stein open set $\Omega\subset X$. 
If $bA\cap K =\emptyset$ and $A\cap K$ is $\cO(A)$-convex
then $A\cup K$ has a fundamental basis of open Stein neighborhoods $\omega$ 
in $X$.
\end{theorem}

%
%
%
%
\begin{figure}[ht]
\psset{unit=0.6cm} 

\begin{pspicture}(-7,-3.5)(7,3.5)
%
%
%
%

\pscircle[linecolor=OrangeRed,fillstyle=crosshatch,hatchcolor=yellow](0,0){3.5}  
\rput(0,3){$\Omega$}						                 

\pscircle[linecolor=DarkBlue,fillstyle=solid,fillcolor=myblue](0,0){2}        
\rput(0,-0.5){$K$}					         	           

\pscurve[linecolor=red,linewidth=1.2pt](-6,0)(-4,0.5)(-3,0)(-2,0)(-1,0.3)(0,0.5)
(1,0.3)(2,-0.2)(3,-0.1)(4,0.2)(6,0)                                           
\psdots[linecolor=red](-6,0)(6,0)

%
%
\psarc[linewidth=0.5pt,linecolor=Black](0,0){2.3}{13}{170}                     
\psecurve[linewidth=0.5pt](2.1,1)(2.22,0.56)(2.8,0.2)(3.4,0.35)(4,0.4)(6,0.25)(6.3,0)(6,-0.2)(4,0)(3.4,-0.2)(3,-0.3)(2.5,-0.4)(2.15,-0.8)(2,-1.5)   						                                             
\psecurve[linewidth=0.5pt](-2,2)(-2.27,0.4)(-2.5,0.2)(-3,0.2)(-3.5,0.6)(-4,0.8)(-6,0.3)(-6.3,-0.05)(-6,-0.2)(-5,0)(-4,0.2)(-3.4,0)(-3,-0.2)(-2.5,-0.3)(-2.2,-0.6)(-2,-1)                                 
\psarc[linewidth=0.5pt](0,0){2.3}{195}{340}                                    

%
%
\rput(-5.5,1.5){$A$}
\psline[linewidth=0.2pt]{->}(-5.5,1.2)(-5.5,0.25)
\rput(4.3,1.5){$\omega$}
\psline[linewidth=0.2pt]{->}(4.3,1.25)(4.3,0.4)
\psline[linewidth=0.2pt]{->}(4,1.4)(2.15,1)
\rput(5.4,1.5){$A$}
\psline[linewidth=0.2pt]{->}(5.4,1.2)(5.4,0.15)

\end{pspicture}
\caption{Theorem \ref{Stein-nbds}} 
\label{Fig1}
\end{figure}

Theorem \ref{Stein-nbds} is the main result of this section
(but see also Theorem \ref{graphs}).
For $X=\C^n$ this follows from  results of 
Wermer \cite{Wer} and Stolzenberg \cite{Stolz}.
We shall only use the special case with $K=\emptyset$, 
but the proof of the general case is not essentially more difficult and 
we include it for future applications.
The necessity of $\cO(A)$-convexity of $K\cap A$ is seen by taking
$X=\C^2$, $A = \{(z,0) \colon |z|\le 3\}$, and 
$K=\{(z,w)\colon 1\le |z|\le 2,\ |w|\le 1\}$:
Every Stein neighborhood of $A\cup K$ 
contains the bidisc $\{(z,w)\colon |z|\le 2,\ |w|\le 1\}$.

In this connection we mention a result of Siu \cite{Siu} 
to the effect that a closed Stein subspace (without boundary) 
of any complex space admits an open Stein neighborhood.
Extensions to the $q$-convex case and  simplifications 
of the proof were given by Col\c toiu \cite{Col0} and Demailly \cite{Demailly}. 
These results do not seem to apply directly to subvarieties with boundaries.

\begin{proof}
We shall adapt the proof of Theorem 2.1 in \cite{FFourier}. 
(It is based on the proof of Siu's theorem \cite{Siu} given in \cite{Demailly}.)
We begin by preliminary results. We have $bA=\cup_{j=1}^m C_j$ 
where each $C_j$ is a closed Jordan curve of class $\cC^2$
(a diffeomorphic image of the circle  $T=\{z\in\C \colon |z|=1\}$).

%
%
%
%
\begin{lemma}
\label{boundary}
There are a Stein open neighborhood $U_j\subset X$ of $C_j$, 
with $\overline U_j\cap K=\emptyset$, and a holomorphic embedding 
$Z=(z,w) \colon U_j\to \C^{1+n_j}$ for some $n_j\in\N$ 
such that $Z(U_j)$ is a closed complex subvariety of the set
\[
	U'_j= \{(z,w)\in \C^{1+n_j} \colon 
	1-r_j <|z| <1+r_j, \  |w_1|<1,\ldots, \ |w_{n_j}|<1 \}
\]
for some $0<r_j<1$, and  
\[
	Z(A\cap U_j)= \{(z,w)\in U'_j \colon z\in \Gamma_j,\  w=g_j(z)\} 
\]
where 
\[
	\Gamma_j =\{ z=re^{i\theta}\in\C \colon 1-r_j < r\le h_j(\theta)\},
\]
$h_j$ is a $\cC^2$ function  close to $1$ (in particular,
$|h_j(\theta)-1|<r_j$ for every $\theta\in\R$), and 
$g_j=(g_{j,1},\ldots,g_{j,n_j}) \colon \Gamma_j\to\disc^{n_j}$ 
is a $\cC^2$ map which is holomorphic in the interior of $\Gamma_j$. 
\end{lemma}

\begin{proof}
We claim that $C_j$, being a totally real submanifold of class $\cC^2$ 
in $X$, admits a basis of open Stein neighborhoods in $X$. 
This is standard when $X$ is smooth (without singularities) 
in which case the squared distance to $C_j$ with respect to any smooth 
Riemannian metric on $X$ is a strongly plurisubharmonic function
in a neighborhood of $C_j$, and its sublevel sets provide a basis
of open Stein neighborhoods of $C_j$. In the general case when $C_j$ contains
some singular points of $X$ we cover $C_j$ by finitely many open 
sets $U_k\subset X$ $(k=1,\ldots,m_j)$ such that each $U_k$ admits
a holomorphic embedding $\phi_k \colon U_k\hra \Omega_k \subset \C^{N_k}$ onto 
a closed complex subvariety $\phi_k(U_k)$ in an open set $\Omega_k \subset \C^{N_k}$.
The function $\rho_k(x)= {\rm dist}^2(\phi_k(x),\phi_k(C_j\cap U_k)) \ge 0$ ($x\in U_j$)
is then strongly plurisubharmonic near the set $\rho_k^{-1}(0) =C_j\cap U_k$.
(We are using the Euclidean distance in the above definition of $\rho_k$.)
Patching these functions $\rho_1,\ldots, \rho_{m_j}$ by a smooth partition of 
unity along $C_j$ in $X$ we obtain a strongly plurisubharmonic function $\rho\ge 0$ 
in a neighborhood of $C_j$ which vanishes precisely on $C_j$,
and the sublevel sets $\{\rho<c\}$ for small $c>0$ provide
a Stein neighborhood basis of $C_j$ \cite{Nar2}. 
The details of the patching argument are similar 
to the nonsingular case and will be omitted.

%
%

Choose a Stein open neighborhood $U_j\Subset X$ of $C_j$.
By shrinking $U_j$ slightly around $C_j$ we may assume that 
$U_j$ embeds holomorphically into a Euclidean space $\C^{1+n_j}$.  
Denote by $C'_j \subset \C^{1+n_j}$ (resp.\ by $A'$) the image of 
$C_j$ (resp.\ of $A\cap U_j$) under this embedding. 
We identify the circle $T$ with $T\times\{0\}^{n_j} \subset \C^{1+n_j}$.
The complexified tangent bundle to $C'_j$, and the complex normal 
bundle to $C'_j$ in $\C^{1+n_j}$, are trivial (since every complex vector bundle
over a circle is trivial). Using standard techniques for totally real 
submanifolds (see e.g.\  \cite{FLO}) we find a $\cC^2$ diffeomorphism 
$\Phi_j$ from a tube around $C'_j$ in $\C^{1+n_j}$
onto a tube around the circle $T$ such that $\Phi_j(C'_j)=T$, and such that
$\dibar\, \Phi_j$ and its total first derivative $D^1(\dibar\, \Phi_j)$ vanish on $C'_j$. 

By Theorems 1.1 and 1.2 in \cite{FLO} we can approximate $\Phi_j$ 
in a tube around $C'_j$ by a biholomorphic map $\Phi'_j$ which maps 
$C'_j$ very close to $T$ and which spreads a collar around $C'_j$ in $A'$ 
as a graph over an annular domain in the first coordinate axis.
Composing the initial embedding $U_j\hra \C^{1+n_j}$ with $\Phi'_j$ 
we obtain (after shrinking $U_j$ around $C_j$) the situation  in the lemma.
\end{proof}

Using the notation in the statement of Lemma \ref{boundary} we set 
\begin{eqnarray}
\label{Lambda-j}
	\Lambda_j &=& \{x\in U_j\colon z(x)\in\Gamma_j\} \subset X, \\
\label{phi-j}
	\phi_j(x) &=& w(x)-g_j(z(x)) \in\C^{n_j}, \quad x\in \Lambda_j.
\end{eqnarray}
We can extend $|\phi_j|^2$ to a $\cC^2$ function on $U_j$ which is
positive on $U_j\bs \Gamma_j$.
Choose additional open sets $U_{m+1},\ldots,U_N$ in $X$ whose closures 
do not intersect any of the sets $U_j\bs \Lambda_j$ for $j=1,\ldots,m$
such that $A\cup K \subset \cup_{j=1}^N U_j$.
By choosing these sets sufficiently small we also get for each 
$j\in \{m+1,\ldots,N\}$ a holomorphic map $\phi_j\colon U_j\to\C^{n_j}$ 
whose components generate the ideal sheaf of $A$ at every point of $U_j$.
If $U_j\cap A=\emptyset$ for some $j$, we take $n_j=1$ and $\phi_j(x)=1$.
Choose slightly smaller open sets $V_j \Subset U_j$ $(j=1,\ldots,N)$
such that $A\cup K \subset \cup_{j=1}^N V_j$. Choose an open set $V\subset X$ with 
$A\cup K\subset V\Subset \cup_{j=1}^N V_j$ and let 
\begin{equation}
\label{Lambda}
	\Lambda = \bigcup_{j=1}^m (\overline V\cap \Lambda_j) 
		\cup \bigcup_{j=m+1}^N (\overline V \cap V_j).
\end{equation}

%
%
%
%
\begin{lemma}
\label{v-delta}
There is a family of $\cC^2$ functions $v_\delta\colon V\to\R$  $(\delta\in (0,1])$
and a constant $M>-\infty$ such that $i\di\dibar \, v_\delta \ge M$ on $\Lambda$
for all $\delta\in (0,1)$, and such that $v_0(x)=\lim_{\delta\to 0} v_\delta(x)$
is of class $\cC^2$ on $V\bs A$ and satisfies $v_0|_A=-\infty$.
\end{lemma}

\begin{proof}
We adapt  the proof of Lemma 5 in \cite{Demailly}.
Let ${\rm rmax}$ denote a regularized maximum 
(p.\ 286 in \cite{Demailly}); this function is increasing and convex 
in all variables (hence it preserves plurisubharmonicity), and it can be 
chosen as close as desired to the usual maximum. 
On every set $V_j$ we choose a smooth function $\tau_j\colon V_j\to\R$ 
which tends to $-\infty$ at $bV_j$. For each $\delta \in [0,1]$  we set
\[
	v_{\delta,j}(x) =  \log(\delta+|\phi_j(x)|^2) + \tau_j(x),\quad x\in V_j,
\]
and $v_\delta(x)={\rm rmax}(\ldots,v_{\delta,j}(x),\ldots)$, where the regularized
maximum is taken over all indices $j\in\{1,\ldots,N\}$ for which $x\in V_j$.
As $\delta\to 0$, $v_\delta$ decreases to $v_0$ and $\{v_0=-\infty\}=A$.
Since the generators $\phi_{j}$ and $\phi_k$ for the ideal sheaf of $A$ 
can be expressed in terms of one another on $U_j\cap U_k$, the quotient
$|\phi_j|/|\phi_k|$ is bounded on $\overline V_j\cap \overline V_k$, and hence 
$(\delta+|\phi_j|^2)/(\delta+ |\phi_k|^2)$ is bounded on $\overline V_j\cap \overline V_k$
uniformly with respect to $\delta\in [0,1]$. Since $\tau_j$ tends to $-\infty$ along
$bV_j$, none of the values $v_{\delta,j}(x)$ for $x$ sufficiently near 
$bV_j$ contributes to the value of $v_\delta(x)$ since the other functions take over in
${\rm rmax}$, and this property is uniform with respect to $\delta \in [0,1]$.
Since $\log(\delta+|\phi_j(x)|^2)$ is plurisubharmonic on  $\Lambda_j$ if 
$j\in \{1,\ldots,m\}$, resp.\ on $U_j$ if $j\in \{m+1,\ldots,N\}$, 
we have $i\di\dibar\, v_{\delta,j} \ge i\di\dibar\, \tau_j$ on the respective sets.
%
%
The above argument therefore gives a uniform lower bound for 
$i\di\dibar\, v_{\delta}$ on the compact set $\Lambda$ (\ref{Lambda}).
However, we cannot control the 
Levi forms of $v_\delta$ from below on the sets $V_j\bs \Lambda_j$ 
for $j\in\{1,\ldots,m\}$ since $\phi_j$ fails to be holomorphic there. 
\end{proof}

%
%
%
%
\begin{lemma}
\label{rho}
Let $U\subset X$ be an open set containing $A\cup K$.
There exists a neighborhood $W$ of $A\cup K$ with 
$\overline W\subset U$ and a $\cC^2$ function $\rho\colon X\to \R$
which is strongly plurisubharmonic on $\overline W$ such that 
$\rho <0$ on $K$ and $\rho >0$ on $bW$.
\end{lemma}

\begin{proof}
Since $A\cap K$ is $\cO(A)$-convex, there exists a compact neighborhood 
$K' \subset U\cap\Omega$ of $K$ such that the set $K'\cap A \subset A\bs bA$ 
is also $\cO(A)$-convex. Since $K$ is  $\cO(\Omega)$-convex, there is a 
smooth strongly plurisubharmonic function $\rho_0\colon \Omega\to\R$ 
such that $\rho_0<0$ on $K$ and $\rho_0>1$ on $\Omega\bs K'$ 
\cite[Theorem 5.1.5, p.\ 117]{Ho}. Set 
$\Omega_c=\{x\in \Omega \colon \rho_0(x)<c\}$. 
Fixing a number $c$ with $0<c< 1/2$ we have  
$K\subset \Omega_c\subset \Omega_{2c}\subset K'$.

Since the restricted function $\rho_0|_{A\cap\Omega}$ is
strongly subharmonic and the set $K'\cap A$ is $\cO(A)$-convex, 
a standard argument \cite[p.\ 737]{FFourier} gives another
smooth function $\wt\rho_0 \colon X\to\R$ which agrees with
$\rho_0$ in a neighborhood of $K'$ in $X$ such that  
$\wt \rho_0|_A$ is strongly subharmonic,
$\wt \rho_0 > c$ on $A\bs \overline \Omega_{c}$,
$\wt \rho_0 > 2c$ on $A\bs \overline \Omega_{2c}$, 
and $\wt\rho_0|_{bA}=c_0\ge 1$ is constant. 

Choose a strongly increasing convex function $h\colon \R\to \R$ 
satisfying $h(t)\ge t$ for all $t\in \R$, $h(t)=t$ for $t\le c$, 
and $h(t)>t+1$ for $t\ge 2c$. The function 
\begin{equation}
\label{rho-1}
	\rho_1 = h\circ \wt \rho_0  \colon X\to\R
\end{equation}
is then strongly plurisubharmonic on $K'$ and along $A$, and it satisfies
\begin{itemize}
\item[(i)]    $\rho_1=\wt \rho_0 =\rho_0$ on $\overline \Omega_c$, 
\item[(ii)]   $\rho_1 \ge \wt\rho_0 >c$ on $A\bs \overline \Omega_{c}$,
\item[(iii)]  $\rho_1 > \wt\rho_0 +1$ on $A\bs \overline \Omega_{2c}$, and 
\item[(iv)]   $\rho_1|_{bA}=c_1 >2$.
\end{itemize}

To complete the proof of Lemma \ref{rho} we shall need the following result;
compare with \cite[Theorem 4]{Demailly}. 

%
%
%
%
\begin{lemma}
\label{psh-extension}
Let $A$ be a compact complex curve with $\cC^2$ boundary in a complex space $X$.
For every function $\rho_1\colon X\to\R$ of class $\cC^2$ such that 
$\rho_1|_A$ is strongly subharmonic there exists a $\cC^2$ function 
$\rho_2\colon X\to\R$ which is strongly plurisubharmonic in a neighborhood 
of $A$ and satisfies $\rho_2|_A=\rho_1|_A$.
\end{lemma}

\begin{proof}
Let $\{U_j\colon j=1,\ldots, N\}$ be the open covering of $A$ 
chosen at the beginning of the proof of Theorem \ref{Stein-nbds}. 
(For the present purpose we delete those sets which do not intersect $A$.)
For each index $j\in\{1,\ldots,m\}$ let 
$Z=(z,w)\colon U_j\to U'_j\subset \C^{1+n_j}$, 
$\Gamma_j$, $\Lambda_j$ and $\phi_j$ be as above.
Denote by $\psi'_j\colon \Gamma_j\times \C^{n_j}\to\R$ the unique
function which is independent of the variable $w\in\C^{n_j}$ 
and satisfies $\rho_1= \psi'_j \circ Z$ on $A\cap U_j$. 
We extend $\psi'_j$ to a $\cC^2$ function $\psi'_j\colon U'_j\to \R$ 
which is independent of the $w$ variable and set 
\begin{equation}
\label{psi-j}
	\psi_j= \psi'_j\circ Z \colon U_j\to \R. 
\end{equation}
Then $\psi_j|_{A\cap U_j}=\rho_1$, and there is an open set 
$\wt \Gamma_j \subset \{1-r_j<|z|<1+r_j\}$, with $\Gamma_j\subset \wt \Gamma_j$, 
such that $\psi_j$ is subharmonic in the open set 
\begin{equation}
\label{wtU-j}
	\wt U_j = \{x\in U_j \colon z(x) \in \wt \Gamma_j \} \subset X.
\end{equation}

By choosing the remaining sets $U_j$ for $j\in \{m+1,\ldots,N\}$
sufficiently small we also get a holomorphic map
$\phi_j\colon U_j\to\C^{n_j}$ whose components generate the ideal
sheaf of $A$ at every point of $U_j$, and a strongly plurisubharmonic 
function $\psi_j\colon U_j\to\R$ extending $\rho_1|_{A\cap U_j}$.

Choose a smooth partition of unity $\{\theta_j\}$ on a neighborhood
of $A$ in $X$ with ${\rm supp}\, \theta_j\subset U_j$ for $j=1,\ldots,N$.
Fix an $\e>0$ and set
\[
	\rho_2(x)= \sum_{j=1}^{N} 
	\theta_j(x)\left( \psi_j(x)+ \e^3 \log\bigl(1+\e^{-4}|\phi_j(x)|^2\bigr) \right).  
\]
For $x\in A$ we have $\rho_2(x)= \sum_j \theta_j(x)\psi_j(x)=\rho_1(x)$.
One can easily verify that $\rho_2$ is strongly plurisubharmonic in a 
neighborhood of $A$ in $X$ provided that $\e>0$ is chosen sufficiently small.
Indeed, as $\e\to 0$, the function
$\e^3 \log\bigl(1+\e^{-4}|\phi_j(x)|^2\bigr)$ is of size $O(\e^3)$,
its first derivative are of size $O(\e)$, and its Levi form at points of 
$A_{reg} \cap U_j$ in the direction normal to $A$ is of size comparable
to $\e^{-1}$, which implies that the Levi form of $\rho_2$ is positive definite
at each point of $A$ provided that $\e>0$ is chosen sufficiently small.
(See the proof of Theorem 4 in \cite{Demailly} for the details.) 
\end{proof}

With $\rho_1$ given by (\ref{rho-1}), and $\rho_2$ furnished by 
Lemma \ref{psh-extension}, we set
\[
	\rho = {\rm rmax} \{\wt\rho_0, \rho_2 -1 \}.
\] 
It is easily verified that $\rho$ is strongly plurisubharmonic on a compact neighborhood 
$\overline W \subset U$ of the set $A\cup \overline \Omega_c$,
$\rho = \wt\rho_0=\rho_0$ on $\overline \Omega_c$ (hence $\rho<0$ on $K$), 
$\rho=\rho_2-1 >\wt\rho_0$ in a neighborhood of $A\bs \Omega_{2c}$,
and $\rho|_{bA}$ has a constant value $C > 1$. 
After shrinking $W$ around $A\cup \overline \Omega_c$ we also have 
$\rho > 0$ on $bW$. 
\end{proof}

%
%
%
%
\smallskip
{\em Completion of the proof of Theorem \ref{Stein-nbds}.}  
We shall use the notation established at the beginning of the proof:
$U_j\subset X$ is an open Stein neighborhood of a boundary curve $C_j \subset bA$,
$\Lambda_j$ and $\phi_j \colon U_j\to\C^{n_j}$ are defined by 
(\ref{Lambda-j}) resp.\ by (\ref{phi-j}), and $\psi_j\colon U_j \to \R$ 
is defined by (\ref{psi-j}).

Let $V$ be an open set containing $A\cup K$, and let $v_\delta\colon V\to\R$ 
$(\delta\in [0,1])$ be a family of functions furnished by Lemma \ref{v-delta}. 
Let $\Lambda$ denote the corresponding set (\ref{Lambda}) 
on which $i\di\dibar\, v_\delta$ is bounded  from below 
uniformly with respect to $\delta\in (0,1]$. As $\delta$ decreases
to $0$, the functions $v_\delta$ decrease monotonically to a function
$v_0$ satisfying $\{v_0=-\infty\}=A$. By subtracting a constant we may assume 
that $v_\delta \le v_1 < 0$ on $K$ for every $\delta \in [0,1]$. 

Given an open set $U\subset X$ containing $A\cup K$, we must 
find a Stein neighborhood $\omega\subset U$ of $A\cup K$.
We may assume that $\overline U\subset V$. 
Let $\rho$ be a function furnished by Lemma \ref{rho}; thus $\rho$ is 
strongly plurisubharmonic on the closure $\overline W\subset U$ 
of an open set $W\supset A\cup K$, $\rho|_K<0$, and
$\rho|_{bW}>0$. Let
\[
	\rho_{\e,\delta} = \rho + \e\, v_\delta\colon \overline W \to\R.
\]
Choose $\e>0$  sufficiently small such that $\rho_{\e,0}>0$ on $bW$
(such $\e$ exists since $\{v_0=-\infty\} =A$); 
hence $\rho_{\e,\delta} \ge \rho_{\e,0} >0$ on $bW$ for every $\delta \in [0,1]$. 
Decreasing $\e>0$ if necessary we may assume that 
$\rho_{\e,\delta}$ is strongly plurisubharmonic on $\Lambda \cap \overline W$
for every $\delta \in (0,1]$ (since the positive Levi form of $\rho$ will compensate 
the small negative part of the Levi form of $\e v_\delta$).  
Fix an $\epsilon$ with these properties. Now choose a sufficiently small 
$\delta>0$ such that $\rho_{\e,\delta} < 0$ on $A$ 
(this is possible since $v_\delta$ decreases to $v_0$ which equals $-\infty$ on $A$). 
Note that $\rho_{\e,\delta} < 0$ on $K$ since both $\rho$ and $v_\delta$ are negative on $K$. 
By continuity $\rho_{\e,\delta}$ is strongly plurisubharmonic
also on the set $\overline W \cap \wt U_j$ for every $j=1,\ldots,m$,
where $\wt U_j \subset U_j$ is an open set of the form  (\ref{wtU-j}).

The function $\psi_j \colon \wt U_j\to\R$ (\ref{psi-j}) is plurisubharmonic on 
the open set $\wt U_j$ (\ref{wtU-j}) which contains $\Lambda_j$,
$\psi_j$ has a constant value $c_1$ on the curve $C_j \subset bA$, and 
$\{\psi_j\le c_1\} = \Lambda_j \supset A\cap U_j$. 
Let $\chi \colon \R\to\R_+$ be a smooth increasing convex function 
with $\chi(t)=0$ for $t\le c_1$ and  $\chi(t)>0$ for $t>c_1$. 
The plurisubharmonic function $\chi\circ \psi_j\colon \wt U_j\to\R$ 
then vanishes on $\Lambda_j$ and is positive on $\wt U_j\bs \Lambda_j$;
extending it by zero along $A$ we obtain a plurisubharmonic function 
$\psi \colon V\to\R_+$ which vanishes of $\overline W\cap \Lambda$ 
and is positive on each of the sets $\wt U_j\bs \Lambda_j$ 
(where it agrees with $\chi\circ \psi_j$). 
By choosing $\chi$ to grow sufficiently fast 
on $\{t> c_1\}$ we can insure that the sublevel set 
\[
	\omega= \{x\in W \colon \psi(x) + \rho_{\e,\delta}(x) <0 \} \Subset W
\]
(which contains $A \cup K$) is contained in the set on 
which $\rho_{\e,\delta}$ is strongly plurisubharmonic.
The purpose of adding $\psi$ is to round off the sublevel set sufficiently
close to $bA$ where it exists from $\Lambda\cap \overline W$, 
thereby insuring that $\omega$  remains in the region where the defining 
function $\psi + \rho_{\e,\delta}$ is strongly plurisubharmonic.
Narasimhan's theorem \cite{Nar2} now implies that $\omega$ is a Stein domain. 
This completes the proof of Theorem  \ref{Stein-nbds}.
\end{proof}

The restriction to one dimensional subvarieties $A\subset X$  
was essential only in the proof of Lemma \ref{boundary}. For higher dimensional
subvarieties we have the following partial result. 

%
%
%
%
\begin{theorem}
\label{graphs}
Let $h\colon X\to S$ be a holomorphic map of a complex space $X$ to a 
complex manifold $S$, and let $D\Subset S$ be a strongly pseudoconvex
Stein domain in $S$. Let $f\colon \bar D\to X$ be a $\cC^2$ section of $h$ 
(i.e., $h(f(z))=z$ for $z\in\bar D$) which is holomorphic in $D$.
If $f(bD)\subset X_{reg}$ and $h$ is a submersion near $f(bD)$
then $A= f(\bar D)$ has a basis of open Stein neighborhoods in $X$.
\end{theorem}

\begin{proof}
The only necessary change in the proof  
is in the construction of the sets $\Lambda_j$ (\ref{Lambda-j}) and
the functions $\phi_j$ (\ref{phi-j}) which describe the subvariety $A\subset X$ 
in a neighborhood of its boundary. When $\dim A = 1$, we could choose 
$\phi_j$ globally around the respective boundary curve $C_j\subset bA$
due to the existence of a Stein neighborhood of $C_j$.  
When $\dim A>1$, this is no longer possible and hence this step 
must be localized as follows.  

Fix a point $p\in bD$ and let $q=f(p) \in bA \subset X_{reg}$. 
Since $h$ is a submersion near $q$, there are local holomorphic coordinates 
$x=(z,w)$ in an open neighborhood $U\subset X$ of $q$, and 
there is an open neighborhood $U'\subset S$ of the point $p=h(q)$ 
such that $h(x)=h(z,w)=z \in U'$ for $x\in U$, 
and $f(z)=(z,g(z))$ for $z\in U'\cap \bar D$. We take 
$\Lambda = \{x=(z,w) \in U\colon z\in U'\cap \bar D\}$
and $\phi(x)=\phi(z,w) = w-g(z)$. Covering $bA$ by finitely many
such neighborhoods, the rest of the proof of Theorem \ref{Stein-nbds}
applies mutatis mutandis.
\end{proof}

%
%
%
%
\begin{corollary}
\label{approximation}
Let $S$ and $X$ be complex manifolds,
and let $D\Subset S$ be a strongly pseudoconvex Stein domain 
with boundary of class $\cC^\ell$. If $2\le r\le \ell$
then every  $\cC^r$ map $f\colon \bar D\to X$ 
which is holomorphic in $D$ is a $\cC^r(\bar D)$ limit
of a sequence of maps $f_j\colon  U_j\to X$ which are holomorphic 
in small open neighborhoods  of $\bar D$ in $S$. 
\end{corollary}

For maps from Riemann surfaces a stronger result is proved in \S 5 below.

\begin{proof} 
When $S=\C^n$, $X=\C^N$, $\ell=2$ and $r=0$, this classical result 
on uniform approximation of holomorphic functions which are continuous
up to the boundary
follows from the Henkin-Ram\'irez integral kernel representation
of functions in $\cA(D)$ (Henkin \cite{Henkin}, Ram\'irez \cite{Ramirez}, 
Kerzman \cite{Kerzman}, Lieb \cite{Lieb}, Henkin and Leiterer \cite[p.\ 87]{HL1}). 
Another  approach which works for $0\le r\le \ell$, $2\le\ell$, is via 
the solution to the $\dibar$-equation with $\cC^r$ estimates 
(Range and Siu \cite{Range-Siu}, Lieb and Range \cite{Lieb-Range0}, 
Michel and Perotti~\cite{Michel-Perotti},
and \cite[Theorem 3.43, VIII/3]{Lieb-Michel}).

Assume now that $X$ is a complex manifold and $2\le r\le \ell$.
By Theorem \ref{graphs} the graph $G_f=\{(z,f(z)) \colon z\in \bar D\}$
admits an open Stein neighborhood $\Omega$ in $S\times X$. 
Choose a proper holomorphic embedding $\psi \colon \Omega\hra\C^N$  and 
a holomorphic retraction $\pi\colon W\to \psi(\Omega)$
from an open neighborhood $W\subset\C^N$ of $\psi(\Omega)$ onto $\psi(\Omega)$. 
Choose a neighborhood $U\subset S$ of $\bar D$ and a sequence of holomorphic maps
$g_j\colon U\to\C^N$ such that the sequence $g_j|_{\bar D}$ 
converges in $\cC^r(\bar D$) to the map $z\to \psi(z,f(z))$ 
as $j\to+\infty$. Denote by $pr_X\colon S\times X \to X$ the projection 
$(z,x)\to x$. Let $U_j=\{z\in U\colon g_j(z)\in W\}$. The sequence
$f_j = pr_X \circ \psi^{-1}\circ \pi \circ g_j \colon U_j\to X$
then satisfies Corollary \ref{approximation}.
\end{proof}

%
%
%
%
\textit{Proof of Theorem \ref{Euclidean-nbds} and Corollary \ref{discs}.}
Let $D\Subset S$ be a smoothly bounded domain in an open
Riemann surface $S$, and let $f\colon\bar D\hra X$ be a $\cC^2$ embedding 
which is holomorphic in $D$. By Theorem \ref{Stein-nbds} the image
$f(\bar D)$ admits an open Stein neighborhood $\Omega\subset X$.
Choose a proper holomorphic embedding $\psi \colon \Omega\hra\C^N$ 
and let $\Sigma=\psi(\Omega)\subset \C^N$.
Also choose a holomorphic retraction $\pi\colon W\to \Sigma$ 
from an open neighborhood $W\subset\C^N$ of $\Sigma$ onto $\Sigma$. 
The embedding $\psi\circ f \colon\bar D\hra \Sigma$ extends
to a $\cC^r$ map $F$ from a neighborhood of $\bar D$ in $S$ to $\Sigma$; 
as $r\ge 2$, $\dibar F$ and its first derivative $D^1(\dibar F)$ vanish on $\bar D$. 

Set $A=F(\bar D) \subset\Sigma$. Let $\nu=T\Sigma|_A/TA$ denote the complex
normal bundle of the embedding $F\colon \bar D\hra \Sigma$; 
this bundle is holomorphic over ${\rm Int}A= F(D)$ and is continuous 
(even of class $\cC^1$) up to the boundary. 
An application of  Theorem B for vector bundles 
which are holomorphic in the interior and continuous up to the boundary
(\cite{Heunemann2}, \cite{Leiterer2}, \cite{Range-Siu})
gives a direct sum splitting $T\Sigma|_A= TA\oplus \nu$ which is 
holomorphic over ${\rm Int}\,A$ and continuous up to the boundary.
(It suffices to follow the proof for vector bundles over open 
Stein manifolds, see e.g.\ \cite[p.\ 256]{GR}.)

Since $A$ is a bordered Riemann surface, the bundle $\nu$ 
is topologically trivial, and hence also holomorphically trivial 
in the sense that it is isomorphic to the product bundle 
$A\times \C^{n-1}$ ($n=\dim X=\dim \Sigma)$ by a
continuous complex vector bundle isomorphism which is holomorphic over the 
interior of $A$  \cite[Theorem 2]{Heunemann1}, \cite{Leiterer1}.
Hence there exist continuous vector fields $v_1,\ldots,v_{n-1}$ tangent
to $\nu \subset T\Sigma|_A$ which are holomorphic in the interior 
of $A$ and generate $\nu$ at every point of $A$.
Considering these fields as maps $A\to T\C^N=\C^N\times\C^N$
we can approximate them uniformly on $A$ by vector fields 
(still denoted $v_1,\ldots, v_{n-1}$) which are holomorphic in 
a neighborhood of $A$ in $\Sigma$ and tangent to $\Sigma$.
(The last condition can be fulfilled by composing them with the 
differential of the retraction $\pi \colon W\to \Omega$.) 
If the approximations are sufficiently close on $A$ then the 
new vector fields are also linearly independent at each point of $A$ and 
transverse to $TA$. The flow $\theta_j^t$ of $v_j$ is defined and holomorphic 
for sufficiently small values of $t\in\C$ beginning at any point near $A$.
The map 
\[
	\wt F(z,t_1,\ldots,t_{n-1})= \theta_1^{t_1}\circ\cdots\circ 
		\theta_{n-1}^{t_{n-1}} \circ F(z)
\]
is a diffeomorphism from an open neighborhood of $\bar D\times\{0\}^{n-1}$
in $S \times \C^{n-1}$ onto an open neighborhood of $A=F(\bar D)$ in 
$\Sigma \subset \C^N$. $\wt F$ is holomorphic in the variables
$t=(t_1,\ldots,t_{n-1})$ and satisfies $\frac {\di \wt F}{\di \bar z}(z,t)=0$ 
for $z\in \bar D$. 

Choose a $\cC^2$ strongly subharmonic function $\rho\colon S\to \R$
such that $D =\{z\in S\colon \rho(z)<0\}$ and $d\rho(z)\ne 0$ for every
$z\in bD=\{\rho=0\}$. For $\e \ge 0$ (small and variable) and 
$M>0$ (large and fixed) the set 
\[
		O_\e = \{(z,t)\in S\times\C^{n-1}\colon \rho(z)+M|t|^2 <\e \}
\]
%
%
is strongly pseudoconvex with $\cC^2$ boundary and is contained 
in the domain of $\wt F$ (the latter condition is achieved by choosing 
$M>0$ sufficiently large). Note that $\bar D\times\{0\}^{n-1} \subset O_\e$ 
for $\e>0$. The properties of $\wt F$ 
described above imply $||\dibar \wt F||_{L^\infty(O_{\e})} = o(\e)$ as $\e\to 0$. 
There are constants $C>0$ and $\e_0>0$ such that for every
$\e \in (0,\e_0)$ the equation $\dibar U=\dibar \wt F$ has a solution 
$U= U_\e \in \cC^2(O_{\e})$ satisfying a uniform estimate 
\begin{equation}
\label{estimate-U}
	||U_\e||_{L^\infty(O_{\e})} \le C||\dibar \wt F||_{L^\infty(O_{\e})}
	= o(\e)
\end{equation}
(see \cite{HL0}, \cite{Lieb-Michel}, \cite{Range-Siu}
and the discussion in \S 3 below). The map 
\[
	G_\e=\pi\circ (\wt F - U_\e) \colon O_\e \to\Sigma \subset\C^N
\] 
is then holomorphic, and it is homotopic to $\wt F|_{O_\e}$ through
the homotopy $G_{\e,s}=\pi\circ (\wt F - sU_\e) \in\Sigma$ $(s\in [0,1])$
satisfying $||G_{\e,s} - \wt F||_{L^\infty(O_{\e})} = o(\e)$
as $\e\to 0$, uniformly in $s\in [0,1]$. 
Choosing $\e>0$ sufficiently small we conclude that 
$G_{\e,s}(z,t)\in \Sigma\bs \wt F(\bar O_0)$  for each  
$(z,t)\in bO_{\e/2}$ and $s\in [0,1]$. 
It follows that for each point $x\in \wt F(\bar O_0)$ 
the number of solutions $(z,t) \in O_{\e/2}$ of the equation 
$G_{\e,s}(z,t) =x$, counted with algebraic multiplicities,
does not depend on $s\in[0,1]$, and hence it
equals one (its value at $s=0$). 
Taking $s=1$ we see that the set  $G_\e(O_{\e/2})$ contains 
$\wt F(\bar O_0) \supset A$. 

From (\ref{estimate-U}) and the interior elliptic regularity 
estimates  \cite[Lemma 3.2]{FLO} we also see that 
$||dU_\e||_{L^\infty(O_{\e/2})} =o(1)$ as $\e\to 0$,
and hence $G_\e$ is an injective immersion on $O_{\e/2}$ 
for every sufficiently small $\e>0$ (since it is a $\cC^1$-small 
perturbation of $\wt F$). For such values of $\e$ 
the set $U_\e:= \psi^{-1}(G_\e(O_{\e/2})) \subset X$ is an open Stein neighborhood 
of $f(\bar D)$, and $U_\e$ is biholomorphic (via $\psi^{-1}\circ G_\e$) to 
the domain $O_{\e/2} \subset S\times\C^{n-1}$.

Since $X$ can be replaced by an arbitrary open neighborhood of $f(\bar D)$
in the above construction, this concludes the proof of Theorem \ref{Euclidean-nbds}.
The same proof gives Corollary \ref{discs}.

%
%
%
%
\section{A Cartan type lemma with estimates up to the boundary}
In this section we prove  one of  our main tools, Theorem \ref{splitting}.

\begin{definition}
\label{Cartan-pair}
A pair of relatively compact open subsets $D_0,D_1 \ss S$ in a complex manifold 
$S$ is said to be a {\em Cartan pair}  of class $\cC^\ell$ $(\ell\ge 2)$ if 
\begin{itemize}
\item[(i)]  the sets $D_0$, $D_1$, $D=D_0\cup D_1$ and $D_{0,1}=D_0\cap D_1$ 
are Stein domains with strongly pseudoconvex boundaries of class $\C^\ell$, and 
\item[(ii)] $\overline {D_0\backslash D_1} \cap \overline {D_1\backslash D_0}=\emptyset$ 
(the separation property). 
\end{itemize}
\end{definition}

Replacing $S$ by a suitably chosen neighborhood of $\overline{D_0\cup D}_1$ 
we can assume that $S$ is a Stein manifold.

Let $P$ be a bounded open set in $\C^n$. 
We shall denote the variable in $S$ by $z$ and the variable in
$\C^n$ by $t=(t_1,\ldots,t_n)$. For each pair of integers $r,s\in\Z_+=\{0,1,2,\ldots\}$ 
we denote by $\cC^{r,s}(\bar D\times P)$ the space of all functions 
$f\colon \bar D\times P\to \C$ with bounded partial derivatives 
up to order $r$ in the $z$ variable and up to order $s$ in the $t$ variable,
endowed with the norm
\[ 
	||f||_{\cC^{r,s}(D\times P)} = 
	\sup\{|D^\mu_z D_t^\nu f(z,t)| \colon z\in \bar D,\ t\in P,\ 
	|\mu| \le r,\ |\nu|\le s\} < +\infty. 
\]	
Here $D^\nu_t$ denotes the partial derivative of order $\nu \in\Z^{2n}$
with respect to the real and imaginary parts of the components 
$t_j$ of $t\in \C^n$. The same definition applies to $D^\mu_z$ 
when $S=\C^m$; in general we cover $\bar D$  by a finite system of local
holomorphic charts $U_j\ss V_j \subset S$, with biholomorphic maps 
$\phi_j\colon V_j\to V'_j \subset \C^m$, and take at each point
$z\in \bar D$ the maximum of the above norms calculated in the $\phi_j$-coordinates
with respect to those charts $(V_j,\phi_j)$ for which $z\in U_j$.
Alternatively, we can measure the $z$-derivatives with respect 
to a smooth Hermitian metric on $S$; the two choices yield 
equivalent norms on $\cC^{r,s}(\bar D\times P)$. Set 
\[
	\cA^{r,s}(D\times P)= \cO(D\times P)\cap \cC^{r,s}(\bar D\times P),
	 \quad r,s\in\Z_+. 
\]
For $t=(t_1,\ldots,t_n)\in\C^n$ we write $|t|=(\sum |t_j|^2 )^{1/2}$.
For a map $f=(f_1,\ldots,f_n)\colon \bar D\times P\to \C^n$ with
components $f_j\in \cC^{r,s}(\bar D\times P)$ we set
\[
		||f||_{\cC^{r,s}(D\times P)} = 
		\biggl( \sum_{j=1}^n ||f_j||^2_{\cC^{r,s}(D\times P)} \biggr)^{1/2} .
\]

Let $\B(t; \delta) \subset \C^n$ denote the ball of radius $\delta>0$ 
centered at $t\in\C^n$. For any subset $P\subset \C^n$ and $\delta>0$ 
we set 
\[
	P_{-\delta}= \{t\in P\colon \B(t; \delta) \subset P\}.
\]

%
%
%
%

\begin{theorem}
\label{splitting}
{\bf (Generalized Cartan's lemma)}
Let $(D_0,D_1)$ be a Cartan pair of class $\cC^\ell$ $(\ell\ge 2)$ 
and let $P$ be a bounded open set in $\C^n$ containing the origin. 
Set $D=D_0\cup D_1$ and $D_{0,1}=D_0\cap D_1$. 
Given $\delta^* >0$ and $r \in\{0,1,\ldots, \ell\}$ there exist numbers 
$\e^* >0$ and $M_{r,s} \ge 1$ $(s=0,1,2,\ldots)$ satisfying the following.
For every map $\gamma\colon \bar D_{0,1} \times P\to\C^n$ of class 
$\cA^{r,0}(D_{0,1} \times P)^n$ satisfying 
\[
	\gamma(z,t)= t + c(z,t), 
	\qquad ||c||_{\cC^{r,0}(D_{0,1} \times P)} < \e^* 
\]
there exist maps $\alpha \colon \bar D_0\times P_{-\delta^*} \to\C^n$,
$\beta \colon \bar D_1 \times P_{-\delta^*} \to\C^n$ of the form 
\[
	\alpha(z,t)=t+a(z,t), \qquad \beta(z,t)=t+b(z,t),
\]
with $a\in \cA^{r,s}(D_0 \times P_{-\delta^*})^n$
and  $b\in \cA^{r,s}(D_1 \times P_{-\delta^*})^n$ for all $s\in \Z_+$,
which are fiberwise injective holomorphic and satisfy   
\begin{equation}
\label{decomposition}
	\gamma(z,\alpha(z,t)) = \beta(z,t), \qquad 
        z \in \bar D_{0,1},\ t \in P_{-\delta^*}
\end{equation}
and also the estimates
\begin{eqnarray*}
	||a||_{\cC^{r,s}(D_0\times P_{-\delta^*})} &\le& 
		M_{r,s} \cdotp ||c||_{\cC^{r,0}(D_{0,1}\times P)}, \\
	||b||_{\cC^{r,s}(D_1\times P_{-\delta^*})} &\le& 
		M_{r,s} \cdotp ||c||_{\cC^{r,0}(D_{0,1}\times P)}.
\end{eqnarray*}
If $\gamma(z,t)=t+c(z,t)$ is tangent to the map $\gamma_0(z,t)=t$
to order $m\in \N$ at $t=0$ (i.e., the function $c(\cdotp,t)$ vanishes to order $m$
at $t=0$) then $\alpha$ and $\beta$ can be chosen
to satisfy the same property.
\end{theorem}

\begin{remark}
The relation (\ref{decomposition}) is equivalent to
\[ 
	\gamma_z = \beta_z \circ \alpha^{-1}_z,\qquad z\in \bar D_{0,1}.
\] 
The classical {\em Cartan's lemma} \cite[p.\ 199, Theorem 7]{GR}
corresponds to the special case when $\alpha_z=\alpha(z,\cdotp)$, $\beta_z$ and 
$\gamma_z$ are linear automorphism of $\C^n$ depending holomorphically on the point 
$z$ in the respective base domain. A version of Cartan's lemma 
without shrinking the base domains was proved by Douady  
\cite{Douady}, and for matrix valued functions
of class $\cA^\infty$ by Sebbar \cite[Theorem 1.4]{Sebbar}. 
Berndtsson and Rosay proved a splitting lemma over the disc $\disc$ 
for bounded holomorphic maps into $GL_n(\C)$ \cite{BR}.
A key difference between all these results and Theorem \ref{splitting}
is that we do not restrict ourselves to fiberwise linear maps.
A result similar to Theorem \ref{splitting}, but less precise 
as it requires shrinking of the base domains,
is  Lemma 2.1 in \cite{ANN} which follows from Theorem 4.1 in \cite{ACTA}.
That lemma does not suffice for the application in this paper
where it is essential that no shrinking be allowed in the base domain.  
\end{remark}

Theorem \ref{splitting} will be proved by a rapidly convergent iteration
similar to the one in the proof of Theorem 4.1 in \cite{ACTA},
but with estimates of derivatives.
At an inductive step we split the map $c(z,t)=\gamma(z,t)-t$ into a difference  $c=b-a$ 
where the maps $a\colon \bar D_0\times P\to\C^n$ and $b\colon \bar D_1\times P\to\C^n$ 
are of class $\cA^{r,0}$, with estimates of their $\cC^{r,0}$ norms in terms of 
the $\cC^{r,0}$ norm of $c$ (Lemma \ref{linear-splitting}).  Set  
\[
	\alpha_z(t)=\alpha(z,t):=t+ a(z,t), \qquad 
	\beta_z(t)=\beta(z,t) :=t+b(z,t).
\]
We then show that for $z\in \bar D_{0,1}$ and $t$ in a smaller set 
$P_{-\delta} \subset\C^n$, with $\epsilon$ sufficiently small compared to $\delta$, 
there exists a map $\wt \gamma \colon \bar  D_{0,1} \times P_{-\delta} \to\C^n$
of the form $\wt\gamma(z,t)=t+\wt c(z,t)$ satisfying
\[
	\gamma_z\circ \alpha_z = \beta_z \circ \wt\gamma_z, \qquad z\in\bar D_{0,1}
\]
and a quadratic estimate
\[ 
	\wt\epsilon = ||\wt c||_{\cC^{r,0}(D_{0,1}\times P_{-\delta})} \le 
	const \cdotp \frac{||c||^2_{\cC^{r,0}(D_{0,1}\times P)} }{\delta}					
\]
(Lemma \ref{inductive-step}).
If $\e=||c||_{\cC^{r,0}(D_{0,1}\times P)}$ is sufficiently small compared to $\delta$
then $\wt \epsilon$ is much smaller than $\epsilon$. 
Choosing a sequence of $\delta\,$'s with the sum $\frac{\delta^*}{2}$ 
and assuming that the initial map $c$ is sufficiently small, the sequences 
of compositions of the maps $\alpha_z$ (resp.\ $\beta_z$), obtained in the individual 
steps, converge on $P_{-\delta^*/2}$ to limit maps $\alpha$ (resp.\ $\beta$)
satisfying $\gamma_z \circ \alpha_z = \beta_z$ for $z \in \bar D_{0,1}$.
After another shrinking of the fiber by  $\frac{\delta^*}{2}$ 
we obtain injective holomorphic maps on $P_{-\delta^*}$ satisfying 
the estimates in Theorem \ref{splitting}.

%
%
%
%
We begin by recalling the relevant results on the solvability of the 
$\dibar$-equation. Let $D$ be a relatively compact strongly pseudoconvex
domain with boundary of class $\cC^\ell$ $(\ell\ge 2)$ in a Stein 
manifold $S$. Let $\cC^{r}_{0,1}(\bar D)$ denote the space of $(0,1)$-forms 
with $\cC^r$ coefficients on $\bar D$, and 
$\cZ^{r}_{0,1}(\bar D) =\{f\in\cC^{r}_{0,1}(\bar D)\colon \dibar f=0\}$.
According to Range and Siu \cite{Range-Siu} and 
Lieb and Range \cite[Theorem 1]{Lieb-Range0} 
(see also \cite[Theorem 1']{Michel-Perotti})
there exists a linear operator $T\colon \cC^0_{0,1}(D)\to \cC^0(D)$ 
satisfying the following properties:
\begin{itemize}
\item[(i)] If $f\in \cC^0_{0,1}(\bar D)\cap \cC^1_{0,1}(D)$ and $\dibar f=0$
then $\dibar (Tf) = f$.
\item[(ii)] If $f\in \cC^0_{0,1}(\bar D) \cap \cC^r_{0,1}(D)$ $(1\le r\le \ell)$  
then for each $l=0,1,\ldots,r$ 
\begin{equation}
\label{T-operator}
	||Tf||_{\cC^{l,1/2}(\bar D)} \le C_l ||f||_{\cC^{l}_{0,1}(\bar D)}.    
\end{equation}
\end{itemize}
The results in \cite{Lieb-Range0} are stated only for the case 
$bD\in \cC^\infty$, but a more careful analysis
shows that one only needs $\cC^\ell$ boundary in order to get 
estimates up to order $\ell$; this is implicitly contained 
in the paper by Michel and Perotti \cite{Michel-Perotti}
(the special case of domains without corners).  
The case of domains in  Stein manifolds easily reduces to the
Euclidean case by standard techniques (holomorphic embeddings 
and retractions). 
Lieb and Range showed that for strongly pseudoconvex domains with smooth
boundaries in $\C^n$ the estimates (\ref{T-operator}) also hold for the
Kohn solution operator $T=\dibar^* N$ (\cite{Lieb-Range1}, \cite[Corollary 2]{Lieb-Range2}). 
Here $\dibar^*$ is the formal adjoint
of $\dibar$ on $(0,1)$-forms (under a suitable choice of a Hermitean metric on $S$)
and $N$ is the corresponding Neumann operator on $(0,1)$-forms on $D$ 
(the inverse of the complex Laplacian $\Box = \dibar \, \dibar^* +\dibar^*\dibar$ 
acting on $(0,1)$-forms).
See also \cite[Theorem 3.43, VIII/3]{Lieb-Michel}; for Sobolev estimates 
see \cite[Theorem 5.2.6, p.\ 103]{Chen-Shaw}.

%
%
%
%
\begin{lemma}
\label{linear-splitting}
Let $D=D_0 \cup D_1\ss S$, $D_{0,1} = D_0\cap D_1$ and $P\subset\C^n$ 
be as in Theorem \ref{splitting}. For every $r\in\{0,1,\ldots,\ell\}$ there are
a constant $C_r\ge 1$, independent of $P$, and linear operators
\[
   A\colon \cA^{r,0}(D_{0,1}\times P)^n \to \cA^{r,0}(D_0\times P)^n, \quad
   B\colon \cA^{r,0}(D_{0,1}\times P)^n \to \cA^{r,0}(D_1\times P)^n
\]
satisfying 
\[
	c= Bc|_{\bar D_{0,1}\times P} - Ac|_{\bar D_{0,1}\times P}, 
	\quad c\in \cA^{r,0}(D_{0,1}\times P)^n,
\] 
and the estimates
\begin{eqnarray*}
       ||Ac||_{\cC^{r,0}(D_0\times P)} &\le& C_r\cdotp ||c||_{\cC^{r,0}(D_{0,1}\times P)}, \\
       ||Bc||_{\cC^{r,0}(D_1\times P)} &\le& C_r\cdotp ||c||_{\cC^{r,0}(D_{0,1}\times P)}.
\end{eqnarray*}
If $c$ vanishes to order $m\in\N$ at $t=0$ then so do $Ac$ and $Bc$. 
\end{lemma}

\begin{proof}
The separation condition (ii) in the definition of a Cartan pair implies that 
there exists a smooth function $\chi$ on $S$ with values in $[0,1]$ such that  
$\chi=0$ in an open neighborhood of $\overline {D_0\backslash D_1}$
and $\chi=1$ in an open neighborhood of $\overline{D_1\bs D_0}$.
Note that $\chi(z)c(z,t)$ extends to a function in $\cC^{r,0}(\bar D_0\times P)$
which vanishes on $\overline {D_0\bs D_1} \times P$, and $\bigl(\chi(z)-1\bigr)c(z,t)$ 
extends to a function in $\cC^{r,0}(\bar D_1\times P)$
which vanishes on $\overline {D_1\bs D_0} \times P$.
Furthermore, $\dibar(\chi c)= \dibar((\chi-1)c)=c\dibar \chi$ is
a $(0,1)$-form on $\bar D$ with $\cC^r$ coefficients and 
with support in $\bar D_{0,1} \times P$, depending
holomorphically on $t\in P$. 

Let $T$ denote a linear 
solution operator to the $\dibar$ equation satisfying (\ref{T-operator}).
For any $c \in \cA^{r,0}(D_{0,1}\times P)$ and $t\in P$ we set
\begin{eqnarray*}
    (A c)(z,t) &=& \bigl(\chi(z)-1\bigr) c(z,t) - T\bigl(c(\cdotp,t) \dibar \chi\bigr) (z), 
    	\quad  z\in\bar D_0. \\   
    (B c)(z,t) &=& \chi(z)c(z,t)  -  T\bigl(c(\cdotp,t) \dibar \chi\bigr) (z), \ \ 
    	\quad\qquad z\in\bar D_1.
\end{eqnarray*}
Then $Ac-Bc=c$ on $\bar D_{0,1}\times P$, $\dibar_z (Ac)=0$, and $\dibar_z(Bc)=0$
on their respective domains. The bounded linear operator  $T$ commutes with 
the derivative $\dibar_t$ on the parameter $t$. 
Since $\dibar_t \left( c(z,t)\dibar \chi(z)\right)=0$, 
we get $\dibar_t (Ac)=0$ and $\dibar_t(Bc)=0$.
The estimates follow from boundedness of $T$ (\ref{T-operator}). 
\end{proof}

%
%
%
%
\begin{lemma}
\label{inductive-step}
Let $D=D_0 \cup D_1\ss S$, $D_{0,1} = D_0\cap D_1$ and $P\subset\C^n$ 
be as in Theorem \ref{splitting}. 
Given $c\in \cA^{r,0}(D_{0,1}\times P)^n$, let $a=Ac$ and 
$b=Bc$ be as in Lemma \ref{linear-splitting}. 
Let $\alpha \colon \bar D_0\times P\to\C^n$,
$\beta \colon \bar D_1\times P\to\C^n$ and 
$\gamma\colon \bar D_{0,1} \times P\to \C^n$ be given by
\[
	\alpha(z,t)=t+a(z,t), \quad \beta(z,t)=t+b(z,t), \quad \gamma(z,t)=t+c(z,t).
\]
Let $C_r \ge 1$ be the constant in Lemma \ref{linear-splitting}. 
There is a constant $K_r>0$ with the following property.
If $4\sqrt{n} C_r ||c||_{\cC^{r,0}(D_{0,1}\times P)} < \delta$ then
there is a map $\wt \gamma \colon \bar  D_{0,1} \times P_{-\delta} \to\C^n$
of the form $\wt\gamma(z,t)=t+\wt c(z,t)$,  with
$\wt c \in \cA^{r,0}(D_{0,1}\times P_{-\delta})^n$, satisfying
the identity
\[
		\gamma_z\circ \alpha_z = \beta_z \circ \wt\gamma_z, \qquad z\in\bar D_{0,1}
\]
and the estimate
\[
	||\wt c||_{\cC^{r,0}(D_{0,1}\times P_{-\delta})} \le 
	K_r \cdotp \frac{||c||_{\cC^{r,0}(D_{0,1}\times P)}^2} {\delta}.
\]
If the functions $a$, $b$ and $c$ vanish to order $m\in\N$ at $t=0$ then so does $\wt c$.  
\end{lemma}

\begin{proof}
We begin by estimating the composition $\gamma_z\circ \alpha_z$.
Since the same estimate will be used below for other compositions
as well, we formulate the result as an independent lemma.

%
%
%
%
\begin{lemma}
\label{estimate-composition}
Let $D$ be a domain with $\cC^1$ boundary in a complex manifold $S$,
let $P$ be an open set in $\C^n$, and let $0<\delta<1$.
Given maps $\alpha_j(z,t)=t + a_j(z,t)$ $(j=0,1)$ with
$a_0 \in \cA^{r,0}(D\times P)^n$, $a_1 \in \cA^{r,0}(D\times P_{-\delta})^n$, 
and $||a_1||_{\cC^{r,0}(D\times P_{-\delta})} < \frac{\delta}{2}$ we have
for all $(z,t)\in \bar D \times P_{-\delta}$
\[
	\alpha_0\bigl(z,\alpha_1(z,t) \bigr)= t + a_0(z,t)+ a_1(z,t) + e(z,t)  
\]
where
\[
	||e||_{\cC^{r,0}(D\times P_{-\delta})} \le 
	\frac{L_r}{\delta}\cdotp ||a_0||_{\cC^{r,0}(D\times P)} 
	\cdotp ||a_1||_{\cC^{r,0}(D\times P_{-\delta})}
\]
for some constant $L_r>0$ depending only on $r$ and $n$.
\end{lemma}

\begin{proof}
We have
\begin{eqnarray*}
	\alpha_0(z,\alpha_1(z,t)) &=& \alpha_1(z,t) + a_0(z,\alpha_1(z,t)) \\
	&=& t + a_1(z,t) + a_0(z,t + a_1(z,t)) \\
	&=& t + a_0(z,t) + a_1(z,t) + e(z,t) 
\end{eqnarray*}
where the error term equals
\[
	e(z,t)= a_0(z,t+a_1(z,t)) - a_0(z,t).
\]
Fix a point $(z,t)\in \bar D  \times P_{-\delta}$.
Since $|a_1(z,t)| < \frac{\delta}{2}$, the line segment
$\lambda\subset \C^n$ with the endpoints $t$ and $\alpha_1(z,t)=t+a_1(z,t)$ 
is contained in $P_{-\delta/2}$. Using the Cauchy estimates 
for the partial derivative $\di_t a_0$ we obtain   
\begin{eqnarray*}
	|e(z,t)| &=& \big| \int_0^1 (\di_t a_0)(z,t+\tau a_1(z,t)) \cdotp a_1(z,t) \, d\tau \big| \\
	&\le&  \sup_{t'\in \lambda} ||\di_t a_0(z,t')|| \cdotp |a_1(z,t)| \\
	&\le&   \frac{2\sqrt{n}}{\delta} \cdotp ||a_0||_{\cC^{0,0}(D\times P)} 
	        \cdotp ||a_1||_{\cC^{0,0}(D\times P_{-\delta})} 
\end{eqnarray*}	
which is the required estimate for $r=0$. We proceed to estimate the 
partial differential of $e(z,t)$:
\begin{eqnarray*}
	\di_z e(z,t) &=& (\di_z a_0)(z,t+a_1(z,t)) - (\di_z a_0)(z,t) + \\
	           && \quad + (\di_t a_0)(z,t+a_1(z,t)) \cdotp (\di_z a_1)(z,t).
\end{eqnarray*}
The difference in the first line equals
\[
	\int_0^1 \di_t (\di_z a_0)(z,t+\tau a_1(z,t)) \cdotp a_1(z,t) \, d\tau
\]
which can be estimated exactly as above (using the Cauchy estimates 
for $\di_t \di_z a_0$)  by 
\[ 
	\frac{const}{\delta} \cdotp ||a_0||_{\cC^{1,0}(D\times P)} 
	        \cdotp ||a_1||_{\cC^{0,0}(D\times P_{-\delta})}. 
\]
Applying the Cauchy estimate for $\di_t a_0$ we estimate the remaining term 
in the expression for $e(z,t)$ by
\[
	\frac{const}{\delta} \cdotp ||a_0||_{\cC^{0,0}(D\times P)} 
	        \cdotp ||a_1||_{\cC^{1,0}(D\times P_{-\delta})}. 
\]
This proves the estimate in  Lemma \ref{estimate-composition} for $r=1$.

We proceed in a similar way to estimate the higher order derivatives of $e$.
In the expression for $\di_z^k e(z,t)$ we shall have a main term 
\[
	(\di_z^k a_0)(z,t+a_1(z,t)) - (\di_z^k a_0)(z,t) =
	\int_0^1 \di_t (\di_z^k a_0)(z,t+\tau a_1(z,t)) \cdotp a_1(z,t) \, d\tau
\]
which is estimated by 
$const\cdotp \delta^{-1} ||a_0||_{\cC^{k,0}(D\times P)} \cdotp 
||a_1||_{\cC^{0,0}(D\times P_{-\delta})}$. 
The remaining terms in $e(z,t)$ are products of partial derivatives
of order $\le k$ of $a_0$ (with respect to both $z$ and $t$ variables) with
partial derivatives of $a_1$ of order $\le k$ with respect to the $z$ variable.
Each $t$-derivative of $a_0$ can be removed by using the Cauchy estimates, 
contributing another $\delta$ in the denominator. The chain rule shows that
each term containing $l$ derivatives of $a_0$ on the $t$ variable 
is multiplied by $l$ factors involving $a_1$ and its $z$-derivatives;
this gives an estimate 
$const\cdotp \delta^{-l} \, ||a_0||_{\cC^{k,0}(D\times P)} \cdotp 
||a_1||^l_{\cC^{k,0}(D\times P_{-\delta})}$. 
Since we have assumed $||a_1||_{\cC^{r,0}(D\times P)} < \frac{\delta}{2}$,
this is less than 
\[
	\frac{const}{\delta} \cdotp ||a_0||_{\cC^{k,0}(D\times P)} \cdotp 
        ||a_1||_{\cC^{k,0}(D\times P_{-\delta})}
\]
and the lemma is proved.
\end{proof}

Let now $\alpha$, $\beta$ and $\gamma$ be as in Lemma \ref{inductive-step}. 
Set $\epsilon= ||c||_{\cC^{r,0}(D_{0,1}\times P)}$;
then $||a||_{\cC^{r,0}(D_0\times P)} \le C_r\e$ and
$||b||_{\cC^{r,0}(D_1\times P)} \le C_r\e$ by  Lemma \ref{linear-splitting}.
Since we have assumed $4\sqrt{n} C_r \e <\delta$, Lemma \ref{estimate-composition} 
with $\alpha_0=\gamma$ and $\alpha_1=\alpha$ gives for 
$z\in \bar D_{0,1}$ and $t\in P_{-\delta}$:
\[
	\gamma(z,\alpha(z,t)) = 
	t+c(z,t)+a(z,t) + e(z,t) = \beta(z,t)+ e(z,t) \in P_{-\delta/2}
\]
where 
\[
	||e||_{\cC^{r,0}(D_{0,1}\times P_{-\delta})} \le 
	\frac{L_r}{\delta}\cdotp ||c||_{\cC^{r,0}(D_{0,1} \times P)} 
	\cdotp ||a||_{\cC^{r,0}(D_{0,1} \times P_{-\delta})}	
	\le \frac{L_r C_r \epsilon^2}{\delta}.
\]

It remains to find a map 
$\wt \gamma(z,t)=t+\wt c(z,t)$ on $\bar D_{0,1}\times P_{-\delta}$
satisfying 
\[
	\beta(z,t)+e(z,t)= \beta(z,t+\wt c(z,t)) = t+\wt c(z,t)+b(z,t+\wt c(z,t))
\]
and an estimate 
\[
	||\wt c||_{\cC^{r,0}(D_{0,1}\times P_{-\delta})} \le const\cdotp \epsilon^2\delta^{-1}.
\]
For the existence of $\wt \gamma$ it suffices to see that the map 
$\beta_z$ is injective on $P_{-\delta/4}$ and $\beta_z(P_{-\delta/4})\supset P_{-\delta/2}$
for every $z\in\bar D_{0,1}$; since $\gamma_z\circ\alpha_z \in P_{-\delta/2}$,
we can then take $\wt\gamma_z=\beta_z^{-1}\circ \gamma_z\circ\alpha_z$.
To see the injectivity of $\beta_z$ note that for $t,t'\in P_{-\delta/4}$, $t\ne t'$, 
we have
\[
	|\beta_z(t)- \beta_z(t')| \ge |t-t'| - |b_z(t)-b_z(t')|
	\ge |t-t'|\left(1 - \frac{4\sqrt{n} C_0\e}{\delta}\right) >0.
\]
(We applied the Cauchy estimate to $\di_t b_z$.)
The inclusion $ P_{-\delta/2} \subset \beta_z(P_{-\delta/4})$
follows from the estimate 
$||b||_{\cC^{r,0}(D_1\times P)} \le C_r\e \le \frac{\delta}{4\sqrt{n}}$
by Rouch\'e's theorem.

In order to estimate $\wt c$ we rewrite its defining equation in the  form
\begin{eqnarray*}
	\wt c(z,t) &=&  b(z,t)-b(z,t+\wt c(z,t)) + e(z,t) \\
	           &=&  -\int_0^1 (\di_t b) \bigl(z,t+\tau \wt c(z,t)\bigr) \cdotp \wt c(z,t) \, d\tau 
	+ e(z,t).
\end{eqnarray*}
Since the path of integration lies in $P_{-\delta/2}$, the Cauchy estimates for $\di_t b$
give
\[
	|\wt c(z,t)| \le \frac{2\sqrt{n} C_0\epsilon}{\delta}\cdotp |\wt c(z,t)| + |e(z,t)| 
	\le \frac{1}{2} \, |\wt c(z,t)| + |e(z,t)| 
\]
and hence $|\wt c(z,t)|\le 2|e(z,t)| \le const\cdotp \e^2\delta^{-1}$.
We proceed inductively to estimate the derivatives $\di_z^k \wt c$ for $k\le r$ 
by differentiating the implicit equation for $\wt c$. 
The top order differential $|\di^k_z \wt c|$ appearing on the right hand 
side is multiplied by a constant $<1$ arising from an estimate
on $b$ (just as was done above); subsuming this term by the left hand side
we obtain the estimates of $|\di^k_z \wt c|$ for all $k\le r$.
Although we obtain a term $\delta^r$ in the denominator, 
we can cancel $r-1$ powers of $\delta$ by appropriate terms 
of size $O(\epsilon)$ just as we did at the end of proof
of Lemma \ref{estimate-composition} to get
$||\wt c||_{\cC^{r,0}(D_{0,1}\times P_{-\delta})} = O(\epsilon^2\delta^{-1})$.  
\end{proof}

%
%
%
%
\noindent 
{\em Proof of Theorem \ref{splitting}.} 
We shall write $(\gamma\alpha)(z,t)= \gamma(z,\alpha(z,t))$,
and similarly for the fiberwise composition of several maps. Let 
\[
	\gamma(z,t)=\gamma_0(z,t)=t+c_0(z,t), \quad
        \epsilon_0 = ||c_0||_{\cC^{r,0}(D_{0,1} \times P)} 
\]
and $\delta^*>0$ be as in Theorem \ref{splitting}. 
We first describe the inductive procedure and subsequently
show convergence provided that $\epsilon_0>0$ is sufficiently small. 
Let $P_0=P$ and $P_*=P_{-\delta^*/2}$. For every $k\in\Z_+$ set 
\[
	\delta_k = 2^{-k-2}\delta^*, \quad P_{k+1}=(P_k)_{-\delta_k}.
\]
Then $\sum_{k=0}^{\infty} \delta_k= \frac{\delta^*}{2}$ and
$\cap_{k=0}^\infty P_k = \bar P_{*}$.
Let $C_r\ge 1$, $K_r\ge 1$ and $L_r \ge 1$ be the constants in
Lemmas \ref{linear-splitting}, \ref{inductive-step} and 
\ref{estimate-composition}, respectively.
We shall inductively construct sequences of maps
\begin{eqnarray*}
	\alpha_k(z,t)    &=& t + a_k(z,t), \quad a_k \in \cA^{r,0}(D_0\times P_k)^n \\ 
        \beta_k(z,t)     &=& t + b_k(z,t), \quad b_k \,\in \cA^{r,0}(D_1\times P_k)^n \\ 
        \gamma_k(z,t)    &=& t + c_k(z,t), \quad c_k \,\in \cA^{r,0}(D_{0,1}\times P_k)^n \\
\end{eqnarray*}
such that, setting 
$
	\epsilon_k = ||c_k||_{\cC^{r,0}(D_{0,1}\times P_k)},
$
the following hold for all $k\in\Z_+$:
\begin{itemize}
\item[($1_k$)] 
$
	||a_k||_{\cC^{r,0}(D_0\times P_k)} \le C_r \epsilon_k, \quad
        ||b_k||_{\cC^{r,0}(D_1\times P_k)} \le C_r \epsilon_k.
$
\item[($2_k$)]  $4\sqrt{n} C_r \epsilon_k < \delta_k = 2^{-k-2}\delta^*$.
\item[($3_k$)]  $\gamma_k\alpha_k=\beta_k\gamma_{k+1}$ on $\bar D_{0,1} \times P_{k+1}$.
\item[($4_k$)]  $\epsilon_{k+1} = ||c_{k+1}||_{\cC^{r,0}(D_{0,1}\times P_{k+1})}  
                \le  K_r \e_k^2 \delta_k^{\,-1} = (4K_r{\delta^{*}}^{-1}) 2^k \epsilon^2_k$.
\end{itemize}
These conditions imply for every $k\in\Z_+$ 
\begin{equation}
\label{together}
	\gamma_0 (\alpha_0\alpha_1\cdots\alpha_{k}) = (\beta_{0}\beta_1\cdots \beta_{k}) \gamma_{k+1}
        \quad {\rm on\ } \bar D_{0,1} \times P_{k+1}. 
\end{equation}
Assuming that $\epsilon_0 = ||c_0||_{\cC^{r,0}(D_{0,1} \times P)} >0$ is sufficiently small
we shall prove that,  as $k\to+\infty$, the sequence of maps
\begin{equation}
\label{tilde-alpha}
	\wt \alpha_k = \alpha_0\alpha_1\cdots\alpha_{k} \colon\bar D_0\times P_k \to\C^n
\end{equation}
converges to a map $\alpha \colon \bar D_0 \times P_*\to\C^n$, the sequence 
\begin{equation}
\label{tilde-beta}
	\wt\beta_k = \beta_{0}\beta_1\cdots \beta_{k}\colon\bar D_1\times P_k \to\C^n
\end{equation}
converges to a map $\beta \colon \bar D_1 \times P_*\to\C^n$, 
and the sequence $\gamma_k$ converges on $\bar D_{0,1} \times P_*$ 
to the map $(z,t)\to t$. (All convergences are in the $\cC^{r,0}$-norms on 
the respective domains.) In the limit we obtain a desired splitting 
\[
	\gamma\alpha=\beta \quad {\rm on} \quad  \bar D_{0,1} \times P_*. 
\]

We begin  at $k=0$ with the given map $\gamma_0(z,t)=t+c_0(z,t)$ on $\bar D_{0,1}\times P_0$. 
Lemma \ref{linear-splitting}, applied to $c_0$, gives maps $a_0$ and $b_0$ satisfying $(1_0)$. 
If ($2_0$) holds (which is the case if $\e_0 = ||c_0||_{\cC^{r,0}(D_{0,1}\times P_0)} >0$ 
is sufficiently small) then Lemma \ref{inductive-step}  furnishes a map 
$\gamma_1\colon \bar D_{0,1} \times P_{1}\to\C^n$ satisfying ($3_0$) and ($4_0$).

Assume inductively that for some $k\in\N$ we already have maps satisfying 
$(1_j)$--$(4_j)$ for $j=0,\ldots,k-1$, and consequently (\ref{together})
holds with $k$ replaced by $k-1$. Lemma \ref{linear-splitting}, applied to 
$c_k(z,t)=\gamma_k(z,t)-t$ on $\bar D_{0,1}\times P_k$, 
gives  maps $a_k$ and $b_k$ satisfying ($1_k$). 
If ($2_k$) holds (and we will show that it does if $\e_0$ is sufficiently small) 
then Lemma \ref{inductive-step}, applied with $\alpha=\alpha_k$, $\beta=\beta_k$, 
$\gamma=\gamma_k$ furnishes a map 
$\wt\gamma= \gamma_{k+1}\colon \bar D_{0,1}\times P_{k+1} \to\C^n$ 
satisfying ($3_k$) and  ($4_k$). This completes the inductive step.

To make the induction work we must insure that the sequence 
$\epsilon_k = ||c_k||_{\cC^{r,0}(D_{0,1}\times P_k)}$ satisfies 
($2_k$) for every $k=0,1,2,\ldots$. To control this process we set 
$N =\max\{ \frac{4K_r}{\delta^{*}}, 1\}$ and define a sequence 
$\sigma_k>0$ by
\begin{equation} 
\label{sigma}
	\sigma_0=\epsilon_0; \quad 	\sigma_{k+1}= 2^k N \sigma_k^2,\ \ k=0,1,2,\ldots
\end{equation}
Any sequence $\e_k \ge 0$ beginning with $\e_0=\sigma_0$ and satisfying ($4_k$) 
for all $k\in\Z_+$ clearly satisfies $\e_k\le \sigma_k$. 
If we can insure (by choosing $\epsilon_0>0$ sufficiently small) that 
\begin{equation} 
\label{estimate1-sigma}
	 \sigma_k < \frac{\delta^*}{2^{k+4}\sqrt{n} C_r},\quad k\in\Z
\end{equation}
then $4\sqrt{n} C_r \e_k \le 4 \sqrt{n} C_r \sigma_k < 2^{-k-2}\delta^* = \delta_k$
and hence ($2_k$) holds.

We look for a solution in the form 
$\sigma_k=2^{\mu_k}N^{\nu_k}{\e_0}^{\tau_k}$. 
From (\ref{sigma}) we get 
\begin{eqnarray*}
 	\mu_{k+1}  &=& 2\mu_k + k, \quad \mu_0=0; \\
	\nu_{k+1}  &=& 2\nu_k + 1, \quad \, \nu_0=0; \\
	\tau_{k+1} &=& 2\tau_k,  \qquad \quad \tau_0=1.
\end{eqnarray*}
Solutions are
\[
	\mu_k=2^k \sum_{l=1}^k l 2^{-l} <2^{k+1}, \quad
	\nu_k = 2^k-1,\quad \tau_k=2^k.
\]
Therefore
\begin{equation}
\label{estimate2-sigma}
	\sigma_k < 2^{2^{k+1}} N^{2^k} {\e_0}^{2^k} =(4N\e_0)^{2^k}, \qquad k\in\N.
\end{equation}
If $\e_0 = ||c_0||_{\cC^{r,0}(D_{0,1}\times P_0)} >0$ is sufficiently small 
then this sequence converges to zero very rapidly and satisfies 
(\ref{estimate1-sigma}). (See Lemma 4.8 in \cite[p.\ 166]{ACTA} for more details.)
For such $\e_0$ we have
\[
	||c_k||_{\cC^{r,0}(D_{0,1}\times P_k)} = \epsilon_k \le \sigma_k \le (4N\e_0)^{2^k} \to 0
\]
and hence $\gamma_k(z,t)\to t$ in $\cC^{r,0}(\bar D_{0,1}\times P_*)$ as $k\to \infty$.

To complete the proof of Theorem \ref{splitting} we must show that the sequences 
(\ref{tilde-alpha}) and (\ref{tilde-beta}) also converge in 
$\cC^{r,0}(\bar D_0\times P_*)$ resp.\ $\cC^{r,0}(\bar D_1\times P_*)$
provided that $\e_0  >0$ is sufficiently small. Write 
\[
	\wt\alpha_k(z,t)= t+\wt a_k(z,t), \quad \wt\beta_k(z,t)= t+\wt b_k(z,t).
\]
By Lemma \ref{estimate-composition} we have
$
	\wt a_{k+1}  = \wt a_k + a_{k+1} + e_{k+1}
$
where 
\[
	||e_{k+1}||_{\cC^{r,0}(D_{0}\times P_{k+1})} \le 
	\frac{L_r}{\delta_k} ||\wt a_k||_{\cC^{r,0}(D_{0}\times P_{k})} 
	||a_{k+1}||_{\cC^{r,0}(D_{0}\times P_{k+1})}. 
\]
Assuming a priori that $||\wt a_k||_{\cC^{r,0}(D_{0}\times P_{k})} \le 1$
for all $k\in\Z_+$  we get the following estimates for the 
$\cC^{r,0}(D_0\times P_{k+1})$ norms:
\[
	||\wt a_{k+1}-\wt a_k||	\le ||a_{k+1}|| + ||e_{k+1}||
	\le C_r\left( 1+ \frac{L_r}{\delta_*} 2^{k+1} \right) \e_{k+1}
	\le R 2^{k+1} \e_{k+1}
\]
with $R=C_r\left(1+ \frac{L_r}{\delta_*} \right)$. 
Note that $\wt a_0=a_0$ and $||a_0||\le C_r \e_0$. Hence 
\[
	||\wt a_0||_{\cC^{r,0}(D_{0}\times P_0)} + 
	\sum_{k=0}^\infty ||\wt a_{k+1} -\wt a_k||_{\cC^{r,0}(D_{0}\times P_{k+1})} 
	\le C_r\e_0 + R\sum_{k=1}^\infty 2^k\e_k. 
\]
Since $\e_k \le \sigma_k \le (4N\e_0)^{2^k}$ for $k\in\N$ (\ref{estimate2-sigma}),
we see that $R \sum_{k=1}^\infty 2^k\e_k < \e_0$ if $\e_0>0$ is sufficiently small.
(See \cite[Lemma 4.8, p.\ 166]{ACTA}  for the details.)
This justifies the assumption $||\wt a_k||_{\cC^{r,0}(D_{0}\times P_{k})} \le 1$
and implies that the sequence $\wt a_k = \wt a_0 + \sum_{j=1}^k (\wt a_j -\wt a_{j-1})$
converges on  $\bar D_0\times P_*$ to a limit $a=\lim_{k\to\infty} \wt a_k$ satisfying 
$||a||_{\cC^{r,0}(D_{0}\times P_{*})} \le (C_0+1)\e_0$. Hence the estimate in 
Theorem \ref{splitting} holds for  $s=0$ with the constant $M_{r,0}=C_0+1$. 

The same proof shows convergence of the sequence 
$\wt b_k \to b$ on $\bar D_1\times P_*$ and the estimate
$||b||_{\cC^{r,0}(D_{1}\times P_{*})} \le (C_0+1)\e_0$.

By shrinking the fiber domain $P_*=P_{-\delta^*/2}$ by an extra
$\frac{\delta^*}{2}$ and applying the Cauchy estimates to the maps
$a(z,\cdotp)$ and $b(z,\cdotp)$ we also obtain
the estimates in the $\cC^{r,s}$ norms in Theorem \ref{splitting}.
In addition, if $\e_0$ is sufficiently small then 
the maps $\alpha(z,\cdotp) \colon P_{-\delta^*} \to\C^n$ 
and $\beta(z,\cdotp) \colon P_{-\delta^*} \to\C^n$
are injective holomorphic for each $z$ in their respective
domain $\bar D_0$ resp.\ $\bar D_1$. 

This completes the proof of Theorem \ref{splitting}.

%
%
%
%
%
\begin{remark}
\label{parametric-case}
Theorem \ref{splitting} holds whenever 
$D_0,D_1, D_{0,1}=D_0\cap D_1, D=D_0\cup D_1$ are relatively compact
domains with $\cC^1$ boundaries satisfying the separation condition  
$\overline {D_0\backslash D_1} \cap \overline {D_1\backslash D_0}=\emptyset$ 
and there exists a linear operator $T\colon \cZ^{r}_{0,1}(\bar D) \to \cC^{r}(\bar D)$
satisfying 
\[
	\dibar (Tf)=f, \qquad ||Tf||_{\cC^{r}(\bar D)} \le C_r ||f||_{\cC^{r}_{0,1}(\bar D)}.
\]
Strong pseudoconvexity of $D_{0,1}$ is not needed here,
but it will be used in the gluing of sprays 
(Proposition \ref{gluing-sprays}). 
The proof of Theorem \ref{splitting} carries over to the 
{\em parametric case} when $\gamma$ depends smoothly 
on real parameters $s=(s_1,\ldots,s_m) \in [0,1]^m \subset\R^m$. 
Indeed, the proof of Lemma \ref{linear-splitting} remains 
valid in the parametric case, and the estimates controlling 
the iteration process are uniform with respect 
to a finite number of $s$-derivatives. This gives a family 
of splittings $\gamma^s_z=\beta^s_z\circ(\alpha^s_z)^{-1}$ for $z\in\bar D_{0,1}$ 
with $\cC^k$ dependence on the parameter $s\in[0,1]^m$ for a given $k\in\N$. 
\end{remark}

%
%
%
%
\section{Gluing sprays on Cartan pairs}
In this section $X$ is an irreducible complex space
and $h\colon X\to S$ is a holomorphic map to a complex manifold $S$. 
Its {\em branching locus} ${\rm br}(h)$  is the union of $X_{sing}$ 
and the set of all those points in $X_{reg}$ at which  $h$ 
fails to be a submersion; thus ${\rm br}(h)$ is an analytic subset of 
$X$, $X'=X\bs {\rm br}(h)$ is a connected complex manifold, 
and $h|_{X'} \colon X'\to S$ is a holomorphic submersion.
For each $x\in X'$ we set $VT_x X =\ker dh_x$, the 
{\em vertical tangent space of $X$}. 

A {\em section} of $h\colon X\to S$ over a subset $D\subset S$ is a map 
$f\colon D\to X$ satisfying $h(f(z))=z$ for all $z\in D$. 
Let $D\Subset S$ be a smoothly bounded domain and $r\in \Z_+$.
A section $f\colon \bar D \to X$ is of class  $\cA^r(D)$
if it is holomorphic in $D$ and $r$ times continuously differentiable
on $\bar D$. (At points of $f(bD) \cap X_{sing}$ we use local holomorphic 
embeddings of $X$ into a Euclidean space.)

%
%
%
%
\begin{definition}
\label{Spray}
{\rm (Notation as above)}
An {\em $h$-spray of class $\cA^r(D)$ with the exceptional set 
$\sigma=\sigma(f)\subset \bar D$ of order $k\ge 0$} is a map $f\colon \bar D\times P\to X$, 
where $P$ (the {\em parameter set} of $f$) is an open subset of 
a Euclidean space $\C^n$ containing the origin, such that the following hold: 
\begin{itemize}
\item[(i)]   $f$ is holomorphic on $D\times P$ and of class $\cC^r$
on $\bar D \times P$, 
\item[(ii)]  $h(f(z,t))=z$ for all $z \in \bar D$ and $t\in P$,
\item[(iii)] 
the maps $f(\cdotp,0)$ and $f(\cdotp,t)$ agree on $\sigma$ up to order $k$ for $t\in P$, and
\item[(iv)]  for every $z\in \bar D\bs \sigma$ and $t\in P$ 
we have $f(z,t)\notin {\rm br}(h)$, and the map
\[
	\di_t f(z,t) \colon T_t \C^n =\C^n  \to VT_{f(z,t)} X
\]
is surjective (the {\em domination condition}).
\end{itemize}
\end{definition}

For a product fibration $h\colon X=S\times Y\to S$, $h(z,y)=z$,
we can identify an $h$-spray $\bar D\times P \to S\times Y$ 
with a {\em spray of maps} $\bar D\times P \to Y$ by 
composing with the projection $S\times Y\to Y$, $(z,y)\to y$.
In this case (ii) is redundant and the domination condition (iv) is replaced by
\begin{itemize}
\item[(iv')] if $z\in \bar D\bs \sigma$ and $t\in P$ then $f(z,t)\in Y_{reg}$ and 
$\di_t f(z,t) \colon T_t \C^n \to T_{f(z,t)} Y$ is surjective.
\end{itemize}

Condition (ii) means that $f_t=f(\cdotp, t)\colon \bar D \to X$ is a 
section of $h$ of class $\cA^r(D)$ for every $t\in P$, 
and by (i) these sections depend holomorphically on the parameter $t$. 
We shall call $f_0$ the {\em core} (or {\em central}\/) {\em section} of the spray. 
Conditions (iii) and (iv) imply that the exceptional set $\sigma(f)$ 
is locally defined by functions of class $\cA^{r}(D)$.

Unlike the sprays used in the Oka-Grauert theory 
which are defined for all values $t\in\C^n$ but are dominant only at the
core section $f_0$, {\em our sprays are local with respect to $t$} 
and dominant at every point $(z,t)$ with $z\notin \sigma$. 
In applications the parameter domain $P$ will be allowed to 
shrink.

%
%
%
%
\begin{lemma}
\label{sprays-exist}
{\bf (Existence of sprays)} 
Let $h\colon X\to S$ be a 
holomorphic map of a complex space $X$ to a complex manifold $S$.
Let $r\ge 2$ and $k\ge 0$ be integers. Let $D$ be a 
relatively compact domain with strongly pseudoconvex boundary 
of class $\cC^2$ in a Stein manifold $S$, and let $\sigma \subset \bar D$ 
be the common zero set of finitely many functions in $\cA^r(D)$. 
Given a section $f_0\colon \bar D\to X$ of class $\cA^r(D)$ 
such that the set $\{z\in\bar D \colon f(z) \in {\rm br}(h)\}$ 
does not intersect $bD$ and is contained in $\sigma$,
there exists an $h$-spray $f\colon \bar D\times P\to X$ of class 
$\cA^r(D)$ with the core section $f_0$ and with the exceptional set 
$\sigma$ of order $k$.
\end{lemma}

\begin{proof}
By Theorem \ref{graphs} there exists a Stein open set $\Omega\subset X$
containing $f_0(\bar D)$. (This is the only place in the proof where
the assumption $r\ge 2$ is used.)
According to \cite[Proposition 2.2]{FORUM} (for manifolds  see \cite[Lemma 5.3]{FP1}) 
there exist an integer $n\in \N$, an open set $V \subset \Omega\times\C^n$ 
containing $\Omega \times \{0\}$, and a holomorphic {\em spray} 
map $s\colon V \to \Omega$ satisfying the following:
\begin{itemize}
\item[(a)]   $s(x,0)=x$ for $x\in \Omega$,
\item[(b)]   $h(s(x,t))=h(x)$ for $(x,t)\in V$,
\item[(c)]   $s(x,t)=x$ when $(x,t)\in V$ and $x\in {\rm br}(h)$, and
\item[(d)]   for each $(x,t)\in V$ with $x\in \Omega\bs {\rm br}(h)$ we have
$s(x,t)\in X\bs {\rm br}(h)$ and the partial differential 
$\di_t s(x,t)|_{t=0} \colon T_0\C^n \to VT_x X =\ker dh_x$
is surjective.
\end{itemize}
A map $s$ with these properties is obtained by composing small complex time 
flows of certain holomorphic vector fields on $\Omega$ which vanish on 
${\rm br}(h)\cap\Omega$ and are tangential to the fibers of $h$.

By the hypothesis we have 
$\sigma=\{z\in \bar D \colon g_1(z)=0,\ldots,g_m(z)=0\}$
where $g_1,\ldots,g_m\in\cA^r(D)$. We can assume that 
$\sup_{z\in \bar D}|g_j(z)| <1$ for $j=1,\ldots,m$.
Denote the coordinates on $(\C^n)^m=\C^{nm}$ by $t=(t_1,\ldots,t_m)$, 
where $t_j = (t_{j,1},\ldots,t_{j,n})\in\C^n$ for $j=1,\ldots,m$.
Let $l\in \N$. The map $\phi_l\colon \bar D \times (\C^n)^m \to \C^n$,
defined by 
\[
	\phi_l(z,t_1,\ldots,t_m)= \sum_{j=1}^m g_j(z)^{k+l}\, t_j,
\]
is a linear submersion $\C^{nm}\to\C^n$ over each point $z\in \bar D \bs\sigma$,
and it vanishes to order $k+l$ on $\sigma$. Let $P \subset \C^{nm}$ be a
bounded open set containing the origin. By choosing the integer
$l$ sufficiently large we can insure that the map 
\[
	f(z,t) = s(f_0(z), \phi_l(z,t)) \in X 
\]
is  a spray $\bar D\times P\to X$ with the core section $f_0$ and 
with the exceptional set $\sigma$ of order $k$. 
All conditions except (iv) are evident.
%
%
To get (iv), let $\Sigma$ denote the set of all points 
$(x,t)\in V$ such that either $x \in {\rm br}(h)$, or $x \notin {\rm br}(h)$
and the map $\di_t s(x,t) \colon T_t\C^n \to VT_{s(x,t)} X$ fails to be 
surjective. Then $\Sigma$ is a closed analytic subset of $V$  
satisfying $\Sigma\cap (\Omega\times\{0\}) = {\rm br}(h)\times\{0\}$
according to property (d) of $s$. (Analyticity of $\Sigma$ 
is clear except perhaps near the points $(x_0,t_0)\in V$ with $x_0\in {\rm br}(h)$.
To see the analyticity near such point we choose a holomorphic 
embedding $\psi\colon U\to \wt U\subset \C^N$ of a small open neighborhood 
$U\subset X$ of $x_0$ onto a local complex subvariety $\wt U=\psi(U) \subset \C^N$
with $\psi(x_0)=0$. Note that $s(x_0,t_0)=x_0$.
There is a holomorphic map $\tilde s$ from a neighborhood of $(0,t_0)\in \C^N\times\C^n$
to $\C^N$ such that $\tilde s(0,t_0)=0$ and $\tilde s(\psi(x),t)=\psi(s(x,t))$;
that is, $\tilde s$ is a local holomorphic extension of $s$ if $U$ is identified
with its image $\wt U\subset \C^N$. Locally near 
the point $(x_0,t_0)$, $\Sigma$ corresponds to the set of points 
$(w,t) \in \C^N\times \C^n$ near $(0,t_0)$ such that $w\in \wt U$ 
and the partial differential $\di_t \tilde s(w,t)$ has rank 
less than $\dim VT(X\bs {\rm br}(h))$; the latter dimension is constant
since $X$ is assumed irreducible. Clearly the latter set is analytic.)
The contact between $\Sigma$ and $\Omega\times\{0\}$ is necessarily of finite 
order along their intersection ${\rm br}(h)\times\{0\}$. By choosing $l\in\Z_+$ large enough 
we insure that $\phi_l(z,t) \in V\bs \Sigma$  for every $z\in\bar D\bs \sigma$ 
and $t\in P$. For such choices $f$ also satisfies the property (iv).
\end{proof}

The following proposition provides the main tool for 
gluing holomorphic sections on Cartan pairs by preserving
their boundary regularity.

%
%
%
%
\begin{proposition}
\label{gluing-sprays}
{\bf (Gluing sprays)} 
Let $h\colon X\to S$ be a holomorphic map from a complex space $X$
onto a Stein manifold $S$. Let $(D_0,D_1)$ be a Cartan pair of
class $\cC^\ell$ $(\ell\ge 2)$ in $S$
(Def.\ \ref{Cartan-pair}) and let $D=D_0\cup D_1$, $D_{0,1}=D_0\cap D_1$. 
Given integers $r\in\{0,1,\ldots, \ell\}$, $k\in\Z_+$, 
and an $h$-spray $f\colon \bar D_0\times P_0\to X$ of class $\cA^r(D_0)$ 
with the exceptional set $\sigma(f)$ of order $k$ and
satisfying $\sigma(f) \cap \bar D_{0,1}=\emptyset$, there is an open 
set $P \ss P_0$ containing $0\in\C^n$ such that the following hold. 

For every $h$-spray $f' \colon \bar D_1 \times P_0 \to X$
of class $\cA^r(D_1)$ with the exceptional set $\sigma(f')$ of order $k$,
with $\sigma(f')\cap \bar D_{0,1}=\emptyset$, such that $f'$ 
is sufficiently $\cC^r$ close to $f$ on $\bar D_{0,1} \times P_0$ 
there exists an $h$-spray $g\colon \bar D \times P \to X$ 
of class $\cA^r(D)$ with the exceptional set $\sigma(g)= \sigma(f)\cup \sigma(f')$
of order $k$ whose restriction $g\colon \bar D_0\times P \to X$ 
is as close as desired to $f \colon \bar D_0\times P \to X$ in the $\cC^r$ topology.
The core section $g_0=g(\cdotp,0)$ is homotopic to $f_0$ on $\bar D_0$, and 
$g_0$ is homotopic to $f'_0$ on $\bar D_1$. In addition,  $g_0$ agrees with 
$f_0$ up to order $k$ on $\sigma(f)$, and $g_0$ agrees with $f'_0$ up to 
order $k$ on $\sigma(f')$.

If $f$ and $f'$ agree to order $m\in\N$ along $\bar D_{0,1}\times\{0\}$ 
then $g$ can be chosen to agree with $f$ to order $m$ along $\bar D_{0}\times\{0\}$,
and to agree with $f'$ to order $m$ along $\bar D_{1}\times\{0\}$.
\end{proposition}

\begin{proof} 
First we find a holomorphic transition map between the two sprays 
(Lemma \ref{transition-map}); decomposing this map by Theorem \ref{splitting} 
we can adjust the two sprays to match them over $\bar D_{0,1}$. 
The first step is accomplished by the following lemma applied
on the strongly pseudoconvex domain $D_{0,1}$.

%
%
%
%
\begin{lemma}
\label{transition-map}
Let $D\ss S$ be a strongly pseudoconvex domain with $\cC^\ell$ boundary $(\ell\ge 2)$
in a Stein manifold $S$, let $P_0$ be a domain in $\C^n$ containing the origin,
and let $f\colon\bar D\times P_0\to X$ be a spray of class 
$\cA^r(D)$ $(0\le r\le \ell)$ with trivial exceptional set. Choose $\e^*>0$. 
There exists an open set $P_1\subset \C^n$, with $0\in P_1\ss P_0$, 
satisfying the following. 
For every spray $f'\colon\bar D\times P_0\to X$ of class $\cA^r(D)$
which approximates $f$ sufficiently closely in the $\cC^r$ topology  
there exists a map $\gamma\colon \bar D\times P_1 \to  \C^n$
of class $\cA^{r,0}(D\times P_1)$ satisfying 
\begin{eqnarray}
\label{gamma1}
	\gamma(z,t) &=& t + c(z,t), \qquad 
		||c||_{\cC^{r,0}(D\times P_1)} < \e^*, \\
\label{gamma2}
         f(z,t) &=& f'(z,\gamma(z,t)),\quad (z,t)\in \bar D \times P_1.
\end{eqnarray}
If $f$ and $f'$ agree to order $m$ along $\bar D\times\{0\}$ then
we can choose $\gamma$ of the form $\gamma(z,t)=t+ \sum_{|J|=m} \wt c_J(z,t)t^J $
with $\wt c_J \in \cA^{r,0}( D\times P_1)^n$.
\end{lemma}

Assuming Lemma \ref{transition-map} for the moment we conclude
the proof of Proposition \ref{gluing-sprays} as follows.
Let $\gamma$ and $P_1$ be as in the conclusion of Lemma \ref{transition-map}
(we emphasize that this lemma is applied on the set $D_{0,1}$).
Choose an open set $P\subset  \C^n$ with $0 \in P\ss P_1$. 
For $\e^*>0$ chosen sufficiently small, Theorem \ref{splitting} 
applied to $\gamma$ gives a decomposition
\begin{equation}
\label{decomp}
	\gamma(z,\alpha(z,t)) = \beta(z,t), \qquad 
        (z,t) \in \bar D_{0,1}\times  P
\end{equation}
where $\alpha \colon \bar D_0\times P \to P_1\subset \C^n$ and 
$\beta \colon \bar D_1 \times P \to P_1\subset \C^n$ are maps of class $\cA^{r,0}$.  
Replacing $t$ by $\alpha(z,t)$ in (\ref{gamma2}) gives 
\begin{equation}
\label{new-spray}
	f(z,\alpha(z,t)) = f'(z,\beta(z,t)), \quad (z,t) \in \bar D_{0,1}\times P.
\end{equation}
Hence the two sides define a map $g\colon \bar D  \times P \to X$ of class 
$\cC^r(\bar D\times P)$ which is holomorphic in $D\times P$. 
Since the maps $\alpha$ and $\beta$ are injective holomorphic 
on the fibers $\{z\} \times P$,  $g$ is a spray with the
exceptional set $\sigma(g)=\sigma(f)\cup \sigma(f')$.

The estimates on $\alpha$ and $\beta$ in Theorem \ref{splitting}
show that their distances from the identity map are controlled by the number
$\epsilon^*$ and hence (in view of Lemma \ref{transition-map})
by the $\cC^r$ distance of $f'$ to $f$ on $\bar D_{0,1}\times P_0$. 
Hence the new spray $g$ approximates $f$ in $\cC^r(\bar D_0\times P)$.
On the other hand, we don't get any obvious control on the $\cC^r$ distance
between $f'$ and $g$ on $\bar D_1\times P$, the problem being that the $\cC^r$ 
norm of $f'$ is not a priori bounded, and precomposing $f'$ by a map 
$\beta$ (even if it is close to the identity map) can still cause a 
big change. However, in our application in \S 6 we shall only need 
to control the range (location) of $g$, and this will be insured by the construction.

Finally, if $f$ and $f'$ agree to order $m$ along $\bar D_{0,1}\times\{0\}$ 
then by Lemma \ref{transition-map} we can choose $\gamma$ of the form 
$\gamma(z,t)=t+ \sum_{|J|=m} \wt c_J(z,t)t^J$ with 
$\wt c_J\in \cA^{r,0}( D_{0,1}\times P_1)^n$ for each multiindex $J$.
Theorem \ref{splitting} then gives a decomposition (\ref{decomp})
where $\alpha(z,t)=t+ \sum_{|J|=m} \wt a_J(z,t)t^J$ 
and $\beta(z,t)=t+\sum_{|J|=m} \wt b_J(z,t)t^J$,
thereby insuring that the spray $g$ (\ref{new-spray}) 
agrees with $f$ resp.\ $f'$ to order $m$ at $t=0$.
This proves Proposition \ref{gluing-sprays}
granted that Lemma \ref{transition-map} holds.

\smallskip
\textit{Proof of Lemma \ref{transition-map}.}
Let $E$ denote the subbundle of $\bar D \times\C^n$ with fibers
\[
	E_z = \ker \left( \di_t f(z,t)|_{t=0} \colon \C^n \to VT_{f(z,0)}X \right),
	\quad z \in \bar D.
\]
This subbundle is holomorphic over $D$ and of class $\cC^r$ on $\bar D$. 
We claim that $E$ is complemented, i.e., there exists a complex vector subbundle 
$G\subset \bar D\times\C^n$ which is continuous on $\bar D$ and holomorphic over $D$ 
such that $\bar D \times\C^n = E\oplus G$. For holomorphic vector bundles on open 
Stein manifolds this follows from Cartan's Theorem B \cite[p.\ 256]{GR};
the same proof applies in the category of holomorphic vector bundles with 
continuous boundary values over a strongly pseudoconvex domain by using the 
corresponding version of Theorem B due to Leiterer \cite{Leiterer2} and 
Heunemann \cite{Heunemann2}. Finally we use a result of 
Heunemann \cite{Heunemann1} to approximate $G$ uniformly 
on $\bar D$ by a holomorphic vector subbundle (still denoted $G$) 
of $U\times \C^n$ over an open neighborhood $U\supset \bar D$; 
a simple proof of this result can be 
found in the Appendix to this paper.
%
%

For each fixed $z\in U$ we write $\C^n\ni t=t'_z\oplus t''_z$
with $t'_z\in E_z$ and $t''_z\in G_z$. The partial differential
$\di_t|_{t=0} f(\cdotp,t)$ gives an isomorphism 
$G|_{\bar D} \to VT_{f_0(\bar D)} X$ and it vanishes on $E$.
The implicit function theorem now gives an open neighborhood $P_1\Subset P_0$ 
of $0\in \C^n$ such that for each spray 
$f'\colon\bar D\times P_0\to X$ which is sufficiently $\cC^r$ 
close to $f$ on $\bar D\times P_0$ there is a unique map 
\[
	\wt \gamma(z,t'_z \oplus t''_z)= t'_z \oplus  (t''_z+ \wt c(z,t)) 
	\in E_z \oplus G_z = \C^n
\]
of class $\cA^{r,0}(D \times P_1)$ solving $f(z,\wt \gamma(z,t))=f'(z,t)$,
and $||\wt c||_{\cA^{r,0}(D_{0,1} \times P_1)}$ is 
controlled by the  $\cC^r$ distance between $f$ and $f'$ on $\bar D\times P_0$.
After shrinking $P_1$ the fiberwise inverse $\gamma(z,t) =t'\oplus (t''_z+ c''(z,t))$ 
of $\gamma$ then satisfies (\ref{gamma2}),
and $||c''||_{\cA^{r,0}(D_{0,1} \times P_1)}$ is controlled by 
the  $\cC^r$ distance between $f$ and $f'$ on $\bar D\times P_0$. 
\end{proof}

\begin{remark}
\label{parametric-gluing}
The additions to Theorem \ref{splitting}, explained in Remark \ref{parametric-case},
yield the corresponding additions to Proposition \ref{gluing-sprays}.
First of all, one can relax the definition of a spray by omitting the
condition regarding the exceptional set. The only essential condition
needed in Proposition \ref{gluing-sprays} is that the spray $f$ is 
{\em dominating on $\bar D_{0,1}$}, in the sense that 
its $t$-differential is surjective on this set at $t=0$.
(This notion of domination agrees with the one introduced by Gromov \cite{Gromov}.)
Approximating such spray $f$ sufficiently closely in the $\cC^r$ topology
on $\bar D_{0} \times P$ (for some open neighborhood $P\subset \C^n$ 
of the origin) by another spray $f'$, we can glue $f$ and $f'$ into
a new spray $g$ over $\bar D_0\cup \bar D_1$ which is dominating
over $\bar D_{0,1}$. The `exceptional set' condition is only needed when
one wishes to interpolate a given spray on a subvariety of $\bar D_0$.
The parametric version of Theorem \ref{splitting} 
(see Remark \ref{parametric-case}) also gives the corresponding
parametric version of Proposition \ref{gluing-sprays} in which the two
$h$-sprays $f$ and $f'$ depend smoothly on a real parameter
$s \in [0,1]^m \subset \R^m$. The remaining ingredients of the proof 
(such as Lemma \ref{transition-map}) carry over to the parametric case 
without difficulties. 
\end{remark}

%
%
%
%

\section{Approximation of holomorphic maps to complex spaces} 
In this section we prove the following approximation theorem for maps
of bordered Riemann surfaces to arbitrary complex spaces.
This result is used in the proof of Theorem \ref{main}
to replace the initial map by another one which maps the boundary
into the regular part of the space.

\begin{theorem}
\label{approximation2}
Let $D$ be a connected, relatively compact, smoothly bounded domain in an 
open Riemann surface $S$, let $X$ be a complex space, 
and let $f\colon \bar D\to X$ be a map of class $\cC^r$ $(r\ge 2)$
which is holomorphic in $D$. Given finitely many points 
$z_1,\ldots,z_l\in D$ and an integer  $k\in\N$, there 
is a sequence of holomorphic maps $f_\nu \colon U_\nu \to X$ in
open sets $U_\nu\subset S$ containing $\bar D$ such that $f_\nu$ agrees 
with $f$ to order $k$ at $z_j$ for $j=1,\ldots, l$ and $\nu\in\N$, 
and the sequence $f_\nu$  converges to $f$ in $\cC^r(\bar D)$ as $\nu\to +\infty$.
If $f(D)\not\subset X_{sing}$, we can also insure that 
$f_\nu(bD)\subset X_{reg}$ for each $\nu\in\N$.
\end{theorem}

\begin{proof}
We proceed by induction on $n=\dim X$. The result trivially holds for $n=0$.
Assume that it holds for all complex spaces of dimension $<n$ for some $n>0$,
and let $\dim X=n$. If $f(D)\subset X_{sing}$ then the conclusion holds by 
applying the inductive hypothesis with the complex space $X_{sing}$. 
Suppose now that $f(D)\not\subset X_{sing}$. The set
\begin{equation}
\label{exc-sigma}
	\sigma = \{z\in \bar D\colon f(z) \in X_{sing} \}
\end{equation}
is compact, $\sigma \cap D$ is discrete, and $\sigma \cap bD$ has empty
relative interior in $bD$. Indeed, as $X_{sing}$ is an analytic
subset of $X$ and hence complete pluripolar, the existence of a nonempty arc 
in $bD$ which $f$ maps to $X_{sing}$ would imply  $f(\bar D)\subset X_{sing}$
in contradiction to our assumption. 

Set $K=\{z_1,\ldots,z_l\}$. Let $bD=\cup_{j=1}^m C_j$ where each $C_j$
is a closed Jordan curve. For each $j=1,\ldots,m$ we choose a point 
$p_j\in C_j \bs \sigma$ and an open set $U_j\subset S$ such that
$p_j\subset U_j$ and $\overline U_j$ does not intersect $\sigma \cup K$.  
We choose the sets $U_j$ so small that $f(\bar D\cap \overline U_j)$
is contained in a local chart of  $X_{reg}$.

\begin{lemma}
\label{f-epsilon}
The map $f$ can be approximated in $\cC^r(\bar D,X)$ by maps
$f'\colon \bar D' \to X$ of class $\cA^r(D',X)$, where $D'\subset S$
is a smoothly bounded domain (depending on $f'$)  satisfying
$D\cup \{p_j\}_{j=1}^m \subset D' \subset D\cup \bigl(\cup_{j=1}^m U_j\bigr)$.
In addition we can choose $f'$ such that it agrees 
with $f$ to order $k$ at $z_j$ for $j\in\{1,\ldots,l\}$. 
\end{lemma}

\begin{proof}
By Theorem \ref{Stein-nbds} the graph of $f$ over $\bar D$ has an open  
Stein neighborhood in $S\times X$. It follows that the set 
$\sigma$ (\ref{exc-sigma}) is the common zero set of finitely many
functions in $\cA^r(D)$. By Lemma \ref{sprays-exist} there is a spray 
$\wt f\colon \bar D\times P\to X$ $(P\subset\C^N)$ 
of class $\cA^r(D)$, with the core map $\wt f(\cdotp,0)=f$ 
and the exceptional set $\tilde \sigma = \sigma\cup K$ of order $k$. 

After shrinking the parameter set 
$P\subset \C^N$ of $\wt f$ around $0\in\C^N$ we may assume that 
$\wt f$ maps the set $E_j=(\overline U_j\cap \bar D)\times \bar P$ 
into a local chart $\Omega\subset X_{reg}$ for each $j=1,\ldots,m$. 
Hence we can approximate the restriction of $\wt f$ to $E_j$
as close as desired in the $\cC^r$ sense by a spray 
$\wt g_j \colon \overline V_j \times P \to X_{reg}$, where 
$V_j$ is an open set in $S$ (depending on $\wt g_j$) 
satisfying $U_j\cap \bar D\subset V_j \subset U_j$. 

If the approximations are sufficiently close, 
Lemma \ref{transition-map} furnishes a transition map 
$\gamma_j$ between $\wt f$ and $\wt g_j$ for each $j$
(we shrink $P$ as needed), and Proposition \ref{gluing-sprays} lets us 
glue $\wt f$  with the sprays $\wt g_j$ into a spray $F$ of class 
$\cA^r(D')$ over a domain $D'\subset S$ as in Lemma \ref{f-epsilon}.
By the construction $F$ approximates $\wt f$ in the $\cC^r(\bar D \times P)$ 
topology, and it agrees with $\wt f$ to a order $k$ at the points $z_j\in K$.
The core map $f'=F(\cdotp, 0) \colon \bar D'\to X$ then satisfies
the conclusion of the lemma. 

A word is in order regarding the application of Proposition \ref{gluing-sprays}.
Unlike in that proposition, the final domain $D'$ in our present situation 
will have to depend on the choices of the sprays $\wt g_j$
(since the size of their $z$-domains in $S$ depends on the rate of approximation).
We can choose from the outset a fixed domain $D_1\subset S$ such that 
$(D, D_1)$ is a Cartan pair in $S$ satisfying
$\overline{D\cap  D}_1 \subset \cup_{j=1}^m (\bar D\cap U_j)$.
Applying Theorem \ref{splitting} gives maps $\alpha$ and $\beta$
over $\bar D$ resp.\ $\bar D_1$; the new spray $F$ is defined 
as $\wt f(z,\alpha(z,t))$ for $z\in \bar D$, and by 
$\wt g_j(z,\beta(z,t))$ for $z\in \bar D_1 \cap U_j$.
Thus we are not using the map $\beta$ on its entire domain of existence, 
but only over the domain of the sprays $\wt g_j$.
\end{proof}

We continue with the proof of Theorem \ref{approximation2}.
Let $f' \colon \bar D' \to X$ be a map furnished by 
Lemma \ref{f-epsilon}. In each boundary curve $C_j\subset bD$ 
we choose a closed arc $\lambda_j\subset C_j$ such that 
$C_j\bs \lambda_j \subset D'$ (this is possible since 
$D'$ contains the point $p_j\in C_j$). 
Let $\xi_j$ be a holomorphic vector field in a neighborhood 
of $\lambda_j$ in $S$ such that $\xi(z)$ points
to the interior of $D$ for every $z\in\lambda_j$.
More precisely, if $D=\{v<0\}$, with $dv\ne 0$ on $bD$,
we ask that $\Re (\xi_j\cdotp v)<0$  on $\lambda_j$;
such fields clearly exist.

Choose a domain $D_0\subset S$ with  $\bar D' \subset D_0$ 
such that $\bar D$ is holomorphically convex in $D_0$. 
(This holds when $D_0\bs \bar D$ is connected.)
The union of $K$ with all the arcs $\lambda_j$ is a 
compact holomorphically convex set in $D_0$. The tangent bundle of $D_0$ is trivial 
which lets us identify vector fields with functions.
Hence there exists a holomorphic vector field $\xi$ on $D_0$
which approximates the field $\xi_j$ sufficiently closely 
on $\lambda_j$ so that it remains inner radial to $D$ there,
and $\xi$ vanishes to order $k$ at the points $z_j\in K$.  
For sufficiently small $t>0$ the flow $\phi_t$ of $\xi$
carries each of the arcs $\lambda_j$ into $D$, 
and hence $\phi_t(\bar D)\subset D'$ provided that $t>0$ is small 
enough. (Recall that $C_j\bs \lambda_j\subset D'$; hence
the points of $\bar D$ which may be carried out of $\bar D$ by the flow
$\phi_t$ along $C_j\bs \lambda_j$ remain in $D'$ for small $t>0$.)

Since the set $\sigma' = \{z\in D' \colon f'(z) \in X_{sing}\}$ 
is discrete, a generic choice of $t>0$ also insures
that $\phi_t(bD)\cap \sigma' =\emptyset$. For such $t$
the map $f'\circ \phi_t$ is holomorphic in an open neighborhood
of $\bar D$, it maps $bD$ to $X_{reg}$, it approximates $f$
in the $\cC^r(\bar D)$ topology, and it agrees with $f$ 
to order $k$ at each point $z_j\in K$. This provides a sequence $f_\nu$ 
satisfying  Theorem  \ref{approximation2}.
\end{proof}

\begin{remark}
\label{Chakrabarti}
D.\ Chakrabarti proved the following approximation result in  \cite[Theorem 1.1.4]{Chak}
(see also \cite{Chak2}):  
{\em If $D$ is a domain in $\C$ bounded by finitely many Jordan curves 
and $X$ is a complex manifold then every continuous map $f\colon \bar D\to X$
which is holomorphic on $D$ can be approximated uniformly on $\bar D$
by maps which are holomorphic in open neighborhoods of $\bar D$ in $\C$.}
A comparison with Theorem \ref{approximation2} shows that there is a 
stronger hypothesis on $X$, but a weaker hypothesis on the map. 
\end{remark}

%
%
%
%

\section{Proof of Theorem \ref{main}}
We begin with the two main lemmas. The induction step in the proof
of Theorem \ref{main} is provided by Lemma \ref{bigstep},
and the key local step is furnished by Lemma \ref{smallstep}.

We denote by $d_{1,2}$ the partial differential with respect to
the first two complex coordinates on $\cn$.

\begin{definition}
\label{2-bump}
Let $A$ and $B$ be relatively compact open sets in a complex space $X$.
We say that {\em $B$ is a 2-convex bump on $A$} (fig.\ \ref{Fig2}) 
if there exist an open set $\Omega\subset X_{reg}$ containing $\bar B$,
a biholomorphic map $\Phi$ from $\Omega$ onto a convex open set
$\omega\subset \C^n$, and smooth real functions $\rho_B\le \rho_A$ on $\omega$ 
such that 
\[
	\Phi(A\cap \Omega)=\{x\in \omega\colon  \rho_A(x)<0\},\ 
	\Phi((A\cup B)\cap \Omega)=\{x\in \omega\colon  \rho_B(x)<0\},
\]
$\rho_A$ and $\rho_B$ are strictly convex with respect to the first 
two complex coordinates, and $d_{1,2}(t\rho_A+(1-t)\rho_B)$ 
is non degenerate on $\omega$ for each $t\in[0,1]$.  
\end{definition}

Let $\rho\colon X\to\R$ be a smooth function which is $(n-1)$-convex 
on an open subset $U\subset X$. If the set 
$\{x\in U\colon c_0\le \rho(x)\le c_1\}$ is compact, 
contained in $X_{reg}$, and it contains no critical points of $\rho$ 
then the set $\{x\in U\colon \rho(x) \le c_1\}$ is obtained 
from $\{x\in U\colon \rho(x) \le c_0\}$ by a finite process
in which every step is an attachment of a 2-convex bump
(Lemma 12.3 in \cite{HL1}). The essential ingredient 
in the proof is Narasim\-han's lemma on local convexification.

The following lemma was proved in \cite{BD3} in the case 
when $X$ is a complex manifold, $D$ is the disc, and for holomorphic maps instead of sprays.
Its proof in \cite{BD3} was based on the solution of the non linear 
Cousin problem in \cite{Ros}. This does not seem to suffice in the case of 
a complex space with singularities and an arbitrary bordered Riemann surface. 
Instead we shall use Proposition \ref{gluing-sprays}. 

Since the complex space $X$ is paracompact, it is metrizable.
Fix a complete distance function $d$ on $X$.

%
%
%
%
%
\begin{lemma}
\label{smallstep}
Let $X$ be an irreducible complex space of $\dim X\ge 2$. 
Let $A\Subset  X$ be relatively compact open subset of $X$ and 
let $B$ be a 2-convex bump on $A$ (Def.\ \ref{2-bump}).
Let $D$ be a bordered Riemann surface with smooth boundary, let $P$ be a 
domain in $\C^N$ containing $0$, and let $k\ge 0$ be an integer.
Assume that $f\colon \bar D\times P\to X$ is a spray of maps
 of class $\cA^2(D)$ with the exceptional set $\sigma$ of order $k$
(Def.\ \ref{Spray}) such that $f_0( bD)\cap \bar A=\emptyset$. 
(Here $f_0=f(\cdotp,0)$ is the core map of the spray.)
Further assume that $K$ is a compact subset of $A$ and 
$U$ is an open subset of $D$ such that $f_0(\bar D\bs U) \cap K=\emptyset$.

Given $\e>0$, there are a domain $P'\subset P$ containing $0\in\C^N$ and 
a spray  of maps $g\colon\bar D\times P'\to X$ of class $\cA^2(D)$, with 
the exceptional set $\sigma$ of order $k$, such that $g_0$ is homotopic 
to $f_0$ and the following hold for all $t\in P'$:
\begin{itemize}
\item[(i)]   $g_t( bD)\cap \overline{A\cup B}=\emptyset$,
\item[(ii)]  $d(g_t(z),f_t( z))<\e$ for $z\in \overline U$,
\item[(iii)] $g_t(\bar D\bs U)\cap K=\emptyset$, and 
\item[(iv)]  the maps $f_0$ and $g_0$ have the same $k$-jets at every point in $\sigma$.
\end{itemize}
\end{lemma}

\begin{proof}
Let  $\Phi\colon X\supset \Omega\to \omega \subset\C^n$ be a biholomorphic map
as in Def.\ \ref{2-bump}. By enlarging the set $U\Subset D$ we may assume 
that $\sigma\subset U$. For small $\l>0$ set
\[
	\omega_\l = \{x\in \omega \colon \rho_B(x)<\lambda, \ \, \rho_A(x)>\lambda\},
		\quad \ \Omega_\lambda = \Phi^{-1}(\omega_\l).
\]
Then $\omega_{\lambda}\ss \omega$ and $\Omega_{\lambda}\ss \Omega$.

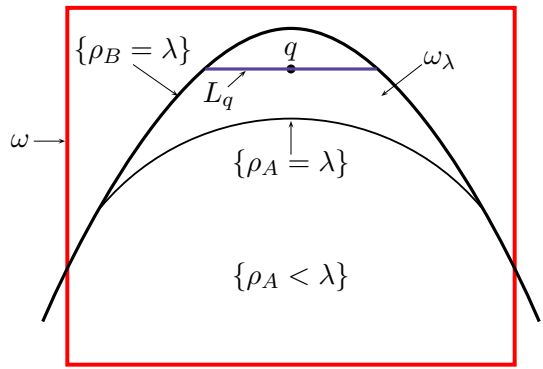
\begin{figure}[ht]
\psset{unit=0.6cm,linewidth=0.7pt}  

\begin{pspicture}(-7,0)(7,8.5)
\psframe[linecolor=red,linewidth=1.5pt](-5,0)(5,8)
\pscurve[linewidth=1.2pt](-5.5,1)(0,7.5)(5.5,1)

\psarc[linewidth=0.7pt](0,0){5.5}{40}{140}
\psdot(0,6.6)
\psline[linewidth=1.2pt,linecolor=Violet](-1.9,6.6)(0,6.6)(1.9,6.6)

\rput(0,2){$\{\rho_A<\lambda\}$}
\rput(-6,5){$\omega$}
\psline[linewidth=0.2pt]{->}(-5.7,5)(-5,5)

\rput(0,7){$q$}
\rput(3.3,6.8){$\omega_\l$}
\psline[linewidth=0.2pt]{<-}(2,6)(2.9,6.7) 
\rput(-1.6,6){$L_q$}
\psline[linewidth=0.2pt]{->}(-1.4,6.2)(-1,6.56) 

\rput(0,4.5){$\{\rho_A=\l\}$}
\psline[linewidth=0.2pt]{->}(0,4.7)(0,5.45)

\rput(-3.5,7){$\{\rho_B=\l\}$}
\psline[linewidth=0.2pt]{->}(-3.5,6.7)(-2.5,6)

\end{pspicture}
\caption{A 2-convex bump}
\label{Fig2}
\end{figure}

Since $f_0( bD)\cap \bar A=\emptyset$, we have $\rho_A(\Phi(f_0(z)))>\lambda$
for every sufficiently small $\l>0$ and for all $z\in bD$ with $f_0(z)\in \Omega$. 
A transversality argument shows that for almost every small
$\l>0$ the set $bD\cap f_0^{-1}(\overline\Omega_{\lambda})$
is a finite union $\cup_{j=1}^{m'} I_j$ of pairwise disjoint closed arcs 
$I_j$ $(j=1,\ldots,m)$ and simple closed curves $I_j$ $(j=m+1,\ldots,m')$.
Fix a $\lambda$ for which the above hold.

If $I_j$ is an arc, we choose a smooth simple closed curve 
$\Gamma_j\subset \bar D\bs U$ such that $\Gamma_j\cap bD$ is a 
\nbd\ of $I_j$ in $bD$, and $\Gamma_j$ bounds a simply connected domain 
$U_j\subset D\bs \overline U$ (fig.\ \ref{Fig3}). Choose  a smooth diffeomorphism 
$h_j\colon\bar\triangle\to\overline U_j$
which is holomorphic on $\triangle$, and choose a compact set
$V_j\subset \overline U_j$ containing a neighborhood of $I_j$ in 
$\bar \disc$.

If $I_j$ is a simple closed curve,  there is a collar neighborhood 
$\overline U_j \subset \bar D\bs \overline U$ of $I_j$ in $\bar D$ 
whose boundary $bU_j=I_j\cup I'_j$ consists of two smooth simple closed curves.
For consistency of notation we set $\Gamma_j=I_j$. 
There are an open subset $W_j$ of $\triangle$  and 
a diffeomorphism $h_j\colon\bar\triangle\bs W_j\to \bar U_j$ which is
holomorphic on $\triangle\bs \overline W_j$ such that $h_j(b\triangle)=\Gamma_j$.
Choose a compact annular neighborhood $V_j$ of $\Gamma_j$ in $U_j\cup \Gamma_j$.

By choosing the sets $U_1,\ldots, U_{m'}$ sufficiently small we can insure that 
their closures are pairwise disjoint and don't intersect $\overline U$, and  
we have
\[
	f_0(\overline U_j)\subset \{x\in \Omega \colon\rho_A(\Phi(x))>\l\},
	\quad j=1,\ldots m'. 
\]

Denote by  $D_1$ the union $\cup_{j=1}^{m'} U_j$.
There is a smoothly bounded open set $D_0$, 
with  $D\bs D_1\subset D_0 \subset D\bs \cup_{j=1}^{m'} V_j$, 
such that $(D_0,D_1)$ is a Cartan pair (Def.\ \ref{Cartan-pair};
see fig.\ \ref{Fig3}).  Let $D_{0,1}=D_0\cap D_1$.

%
%
%
%
\begin{figure}[ht]
\psset{unit=0.6cm, linewidth=0.7pt}  
\begin{pspicture}(-6,-6)(6,6.2)

\pscircle[linewidth=0.8pt](0,0){5}

%
%
%
%
\pscustom[linestyle=none,fillstyle=crosshatch,hatchcolor=yellow,linewidth=0pt]
{
\psarc(0,0){5}{45}{135}
\psarc(-3.18,3.18){0.5}{135}{315}
\psarc(0,0){4}{45}{135}
\psarc(3.18,3.18){0.5}{-135}{45}
}

%
%
\psarc[linecolor=Violet,linewidth=1pt](0,0){5}{45}{135}
\psarc[linecolor=Violet,linewidth=1pt](0,0){4}{45}{135}
\psarc[linecolor=Violet,linewidth=1pt](-3.18,3.18){0.5}{135}{315}
\psarc[linecolor=Violet,linewidth=1pt](3.18,3.18){0.5}{-135}{45}

%
%
\psarc[linecolor=black,linewidth=1.4pt,arrows=*-*](0,0){5}{70}{110}
\psecurve(-4,3)(-3,4)(0,4.5)(3,4)(4,3)

%
%
\pscustom[linestyle=none,fillstyle=crosshatch,hatchcolor=yellow]
{
\psarc(0,0){5}{165}{255}
\psarc(-1.165,-4.347){0.5}{255}{75}
\psarc(0,0){4}{165}{255}
\psarc(-4.347,1.165){0.5}{-15}{165}
}

\psarc[linecolor=Violet,linewidth=1.2pt](0,0){5}{165}{255}
\psarc[linecolor=Violet,linewidth=1.2pt](0,0){4}{165}{255}
\psarc[linecolor=Violet,linewidth=1.2pt](-4.347,1.165){0.5}{-15}{165}
\psarc[linecolor=Violet,linewidth=1.2pt](-1.165,-4.347){0.5}{255}{75}

\psarc[linecolor=black,linewidth=1.4pt,arrows=*-*](0,0){5}{190}{230}
\psecurve(4.598,-1.96)(-1.964,-4.5981)(-3.897,-2.25)(-4.964,0.598)(-4.598,1.964)

%
%
\pscustom[linestyle=none,fillstyle=crosshatch,hatchcolor=yellow]
{
\psarc(0,0){5}{285}{15}
\psarc(4.347,1.165){0.5}{15}{195}
\psarc(0,0){4}{285}{15}
\psarc(1.165,-4.347){0.5}{105}{285}
}

\psarc[linecolor=Violet,linewidth=1pt](0,0){5}{285}{15}
\psarc[linecolor=Violet,linewidth=1pt](0,0){4}{285}{15}
\psarc[linecolor=Violet,linewidth=1pt](4.347,1.165){0.5}{15}{195}
\psarc[linecolor=Violet,linewidth=1pt](1.165,-4.347){0.5}{105}{285}

\psarc[linecolor=black,linewidth=1.4pt,arrows=*-*](0,0){5}{310}{350}
\psecurve(-4.598,-1.96)(1.964,-4.5981)(3.897,-2.25)(4.964,0.598)(4.598,1.964)

%
%
\pscircle[linestyle=none,fillstyle=solid,fillcolor=white](0,0){4}

%
%
\psline[linewidth=0.2pt]{->}(-2.3,2.8)(-2.8,3.3) 
\rput(-1.9,2.4){$U_1$}

\psline[linewidth=0.2pt]{<-}(-4.3,0)(-3.4,0) 
\rput(-2.9,0){$U_2$}

\psline[linewidth=0.2pt]{->}(3.4,0)(4.3,0)  
\rput(2.9,0){$U_3$}

\rput(0,0){$D_0$}

\rput(0.1,2.9){$\Gamma_1$}
\psline[linewidth=0.2pt]{->}(0,3.35)(0,4) 

\rput(-2.5,-1.4){$\Gamma_2$}
\psline[linewidth=0.2pt]{->}(-2.95,-1.6)(-3.5,-2)

\rput(2.5,-1.4){$\Gamma_3$}
\psline[linewidth=0.2pt]{->}(2.95,-1.6)(3.5,-2)

\rput(0,5.9){$I_1$}
\psline[linewidth=0.2pt]{->}(0,5.5)(0,5.05)

\rput(-5.1,-2.85){$I_2$}
\psline[linewidth=0.2pt]{->}(-4.9,-2.7)(-4.4,-2.4)

\rput(5.15,-2.85){$I_3$}
\psline[linewidth=0.2pt]{->}(4.9,-2.7)(4.4,-2.4)

\rput(-5,3.2){$bD$}
\psline[linewidth=0.2pt]{->}(-4.75,3.1)(-4.2,2.8)

\rput(0,-2.5){$bD_0$}
\psline[linewidth=0.2pt]{->}(-0.6,-2.5)(-3.5,-3)
\psline[linewidth=0.2pt]{->}(0.55,-2.5)(3.5,-3)
\psline[linewidth=0.2pt]{->}(0,-2.8)(0,-5)

\rput(-5,5){$D_1=\cup_j U_j$}

\end{pspicture}
\caption{Cartan pair $(D_0,D_1)$}
\label{Fig3}
\end{figure}
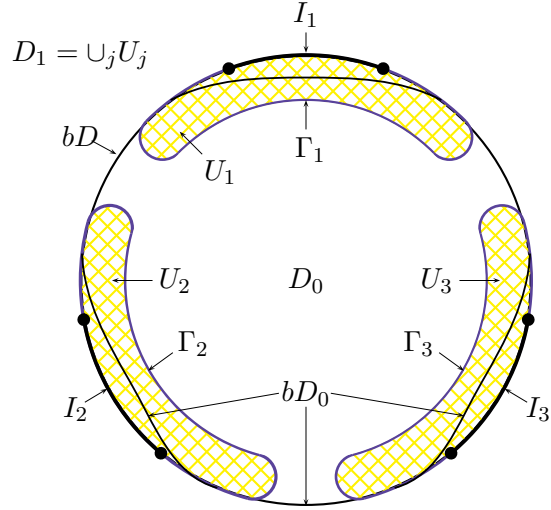

Our goal is to approximate  $f$ in the $\cC^2$ topology on 
$\bar D_{0,1}$ by a spray $f'$ over $\bar D_1$ such that 
the maps $f'_t$ will satisfy properties (i) and (iii) on its domain.
(The final spray $g$ over $\bar D$ will be obtained by gluing 
the restriction of $f$ to $\bar D_0$ with the spray $f'$, 
using Proposition \ref{gluing-sprays}.)
To this end we shall now find a suitable family of
holomorphic discs which will be used to increase 
the value of $\rho\circ f_0$ on the part of $bD$ which is mapped
by $f_0$ into $\Omega_\l$.

Consider the homotopy $\rho_s\colon \omega\to \R$
defined by 
\[
	\rho_s=(1-s)(\rho_A-\lambda)+s(\rho_B-\lambda),\quad s\in[0,1].
\]
The function $\rho_s$ is strictly convex with respect to the first 
two coordinates (since it is a convex combination of functions with this property), 
and $d_{1,2}\rho_s$ is non degenerate on $\omega$ by the definition of
a 2-convex bump. As the parameter $s$ increases from $s=0$ to $s=1$, the  sets
$\{\rho_s \le0\}$ increase smoothly from $\{\rho_A\le \l\}$ to $\{\rho_B\le\l\}$.
(Inside $\omega_\l$ these sets are strictly increasing.)
For each point $q\in \omega_\l$ we have
$\rho_A(q) > \lambda$ while $\rho_B(q) <\lambda$;
hence there is a unique $s\in [0,1]$ such that $\rho_s(q)=0$.
Write $q =(q_1,q_2,q'')$, with $q''\in\C^{n-2}$. The set 
\[
	M_{s,q''} = \{(x_1,x_2,q'')\in \omega \colon \rho_s(x_1,x_2,q'')=0\}
\]
is a real three dimensional submanifold of $\C^2\times\{q''\}$.
Let $T_q M_{s,q''}$ denote its real tangent space at $q$;
then $E_q=T_q M_{s,q''}\cap i\, T_q M_{s,q''}$ is a complex line
in $T_q \C^n=\C^n$. By strict convexity of $\rho_B$ with respect to the 
first two variables the intersection 
\[
	L_q= (q+E_q) \cap \{x\in \omega \colon \rho_B(x)\le \l\}
\]
is a compact, connected, smoothly bounded convex subset of $q+E_q$
with $bL_q\subset \{\rho_B=\l\}$ (fig.\ \ref{Fig2}). 
The sets $L_q$ depend smoothly on $q\in\omega_\l$ and 
degenerate to the point $L_q=\{q\}$ for
$q\in b\,\omega_\l \cap\{\rho_A>\l\}$. We set $L_q=\{q\}$ 
for all points $q\in \omega$ with $\rho_B(q)\ge \l$.

Given a point $z\in \Gamma_j \subset bD_1$ for some $j\in\{1,\ldots,m'\}$,
we set
\[
	\wt L_z = L_q \ \ \text{with}\ \  q=\Phi(f_0(z)). 
\]
The definition is good since $\rho_A(\Phi(f_0(z)))>\l$ for all $z\in \bar D_1$.

An elementary argument (see e.g.\ \cite[Section 4]{Glo}) gives for each 
$j\in\{1,\ldots m'\}$ a continuous map $H_j\colon \Gamma_j \times\bar\triangle\to \omega$
such that for each $z\in I_j$ the map $\bar\disc\ni \eta \mapsto H_j(z,\eta) \in \wt L_z$
is a holomorphic parametrization of $\wt L_z$ and $H_j(z,0) = \Phi(f_0(z))$; 
if $z\in\Gamma_j\bs I_j$ then $H_j(z,\eta) = \Phi(f_0(z))$ for all $\eta\in\bar\disc$. 

Recall that $h_j$ is a parametrization of $\overline U_j$ by a 
$\bar \disc$ if $j\in\{1,\ldots,m\}$, resp.\ by an annular region in $\bar \disc$ 
if $j\in\{m+1,\ldots,m'\}$. Let 
$G_j\colon b\triangle\times\bar\triangle\to\C^n$ be defined by 
\[
	G_j(\zeta,\eta)=H_j(h_j(\zeta),\eta)-\Phi(f_0(h_j(\zeta))), 
	\quad \zeta \in b\disc,\ \eta\in\bar\disc.
\]
Observe that $G_j(\zeta,\eta)=0$ if $\zeta \in h_j^{-1}(\Gamma_j\bs I_j)$ 
and $\eta\in\bar\disc$.

Let $\B \subset \C^n$ denote the unit ball and $\d\,\B$ the ball of radius $\d$. 
For each $j\in\{1,\ldots,m'\}$ and each $\d>0$ we solve approximately the 
Riemann-Hilbert problem for the map $G_j$, using \cite[Lemma 5.1]{Glo}, to
obtain a holomorphic polynomial map $Q_{\delta,j}\colon \C\to\C^n$ 
satisfying the following properties:
\begin{eqnarray}
\label{P1.1}
	 Q_{\delta,j}(\zeta)  &\in&  G_j(\zeta,b\triangle)+\d\,\B \ \ \text{ for } \zeta \in b\triangle,\\
\label{P1.2} 
	|D^2 Q_{\delta,j}(\zeta)| &<& \d \ \ \text{ for } \zeta\in h_j^{-1}(\overline {U_j\bs V_j}), \\
\label{P1.3} 
	 Q_{\delta,j}(\zeta)  &\in&  G_j(b\triangle,\bar\triangle)+ \d\,\B \ \ 
	 	\text{ for } \zeta \in h_j^{-1}(\overline U_j).
\end{eqnarray}
Here $D^2 Q=(Q,Q',Q'')$ is the second order jet of $Q$. 
Although Lemma 5.1 in \cite{Glo} only gives a uniform estimate in (\ref{P1.2}),
we can apply it to a larger disc containing $h_j^{-1}(\overline {U_j\bs V_j})$ 
in its interior to obtain the estimates of derivatives.

Define a map $Q_{\delta}\colon \bar D_1=\cup_{j=1}^{m'} \overline U_j \to\C^n$ by 
\[
		Q_{\delta}(z) = Q_{\delta,j}(h_j^{-1}(z)), \quad z\in \overline U_j. 
\]
By (\ref{P1.2}) the map $Q_{\delta}$ and its first two derivatives have modulus bounded by $\d$ 
on $\cup_{j=1}^{m'} \overline {U_j\bs V_j}$, and hence on $\bar D_{0,1}$.  
If $z\in \Gamma_j\cap bD$ then (\ref{P1.1}) gives  
\[
		|Q_\delta(z)+ \Phi(f_0(z)) - H_j(z,\eta)| < \delta \ \ \text{for some}\ \ \eta\in b\disc,
\]
and hence the point $Q_\delta(z)+ \Phi(f_0(z))$ is contained in the $\d$-neighborhood of 
$b\wt L_z$. Recall that for $z\in I_j$ we have 
$b\wt L_z \subset \{\rho_B=\lambda\}$, and for $z\in \Gamma_j\bs I_j$
we have $\wt L_z= \{ \Phi(f_0(z)) \}$. By choosing $\d_0>0$ sufficiently small 
we insure that  
\[
		\rho_B\bigl( Q_{\delta}(z)+ \Phi(f(z,t)) \bigr) > 0
\]
for all $z\in \Gamma_j\cap bD$, $j=1,\ldots,m'$, $0< \d<\d_0$, and all $t$ in a certain 
neighborhood $P_0\subset P$ of $0\in\C^N$. 
For such choices (and a fixed $\delta\in(0,\delta_0)$) 
the map $f'=f'_\delta \colon \bar D_1\times P_0\to X$, defined by
\[
		f'(z,t) = \Phi^{-1}\bigl( Q_\delta(z) + \Phi(f(z,t)) \bigr),
		\quad z\in \bar D_1,\ t\in P_0,
\]
is a spray of maps of class $\cA^2(D_1)$, with trivial (empty) exceptional set,
whose  boundary values on $bD_1\cap bD$ lie outside of $\overline {A\cup B}$.
By choosing $\d>0$ small enough we insure that $f'$ approximates the spray $f$ 
as close as desired in the $\cC^2$ norm on $\bar D_{0,1} \times P_0$. 

By Proposition \ref{gluing-sprays} we can glue $f$ and $f'$ into a spray of maps 
$g \colon \bar D\times P'\to X$ approximating  $f$ on $\bar D_0\times P'$; 
hence the central map $g_0=g(\cdotp,0)$ satisfies property (ii) in Lemma \ref{smallstep},
and also property (i) on $bD_0\cap bD$.
For $z\in \bar D_1$ we have $g(z,t)=f'(z,\beta(z,t))$ by (\ref{new-spray}),
where the $\cC^2$ norm of $\beta$ is controlled by $\d$.
Choosing $\d >0$ sufficiently small we insure that for each $z\in bD_1\cap bD$ 
we have $g_0(z)=g(z,0) \in X\bs \overline {A\cup B}$, so (i) holds also on $bD_1\cap bD$.
Similarly, since $f'_t(\bar D_1)$ does not intersect $\bar A \supset K$,
we see that $g_0$ satisfies property (iii). 
By shrinking $P'$ we obtain the same properties for all
maps $g_t$, $t\in P'$. Finally, property (iv) holds by the construction 
(this does not depend on the choice of the constants).
\end{proof}

%
%
%
%

\begin{lemma}
\label{bigstep}
Let $X$ be an irreducible complex space of dimension $n\ge 2$, and let $\rho\colon X\to\R$ 
be a smooth exhaustion function 
which is $(n-1)$-convex on $\{x\in X\colon\rho(x)>M_1\}$.
Let $D$ be a finite Riemann surface, let $P$ be an open set in $\C^N$ containing $0$,
and let $M_2>M_1$. Assume that $f\colon\bar D\times P\to X$ is a spray of maps 
of class $\cA^2(D)$ with the exceptional set $\sigma\subset D$ of order $k\in\Z_+$, 
and $U\ss D$ is an open subset such that 
$f_0(z)\in \{x\in X_{reg}\colon\rho(x)\in (M_1,M_2)\}$ for all $z\in\bar D\bs U$. 
Given $\e>0$ and a number $M_3>M_2$, there exist a domain $P'\subset P$ containing $0\in\C^N$ and
a spray of maps $g\colon\bar D\times P'\to X$ of class $\cA^2(D)$, 
with exceptional set $\sigma$ of order $k$, satisfying the following properties:
\begin{itemize}
\item[(i)]   $g_0(z)\in \{x\in X_{reg}\colon\rho(x)\in (M_2,M_3)\}$ for $ z\in bD$,
\item[(ii)]  $g_0(z)\in \{x\in X\colon\rho(x)>M_1\}$ for $z\in\bar D\bs U$,
\item[(iii)] $d(g_0(z),f_0(z))<\e$ for $z\in \overline U$, and
\item[(iv)]  $f_0$ and $g_0$ have the same $k$-jets at each of the points in $\sigma$.
\end{itemize}
Moreover, $g_0$ can be chosen homotopic to $f_0$. 
\end{lemma}

\begin{proof}
The idea is the following. Lemma \ref{smallstep} allows us to push the boundary of 
our curve out of a 2-convex bump in $X$. By choosing these bumps carefully
we can insure that in finitely many steps we push the boundary of the curve 
to a given higher super level set of $\rho$ (property (i)); at the same time we take care
not to drop it substantially lower with respect to $\rho$ (property (ii)) and to approximate the given map on the compact subset $\overline U \subset D$ (property (iii)). 
In the construction we always keep the 
boundary of the image curve in the regular part of $X$. Special care must
be taken to avoid the critical points of $\rho$.
We now turn to details. 

By \cite[Lemma 5]{Demailly} there exists an {\em almost \psh\ function} $v$ on $X$ 
(i.e., a function whose Levi form has bounded negative part on each compact in $X$)
which is smooth on $X_{reg}$ and satisfies $v=-\infty$ on $X_{sing}$. 
We may assume that $v<0$ on $\{\rho \le M_3+1\}$. 

For every sufficiently small $\d>0$ the function 
$\tau_\d=\rho-M_1+ \d v$ is $(n-1)$-convex on $\{\rho \le M_3\}$,
and its Levi form is positive on the linear span of the eigenspaces 
corresponding to 
the positive eigenvalues of the Levi form of $\rho$ at each point. 
Note that $X_{sing} \cup \{\rho \le M_1\} \subset \{\tau_\d<0\}$.
Since $\rho(f_0(z))>M_1$ and $f_0(z)\in X_{reg}$ for all $z\in bD$, 
we have $\tau_\d(f_0(z))>0$ for all $z\in bD$ and all small $\d>0$.
Fix $\d>0$ for which all of the above hold and write $\tau=\tau_\d$.

Choose a number $M\in (M_2,M_3)$. (The central map $g_0$ 
of the final spray will map $bD$ close to 
$\{\rho=M,\ \tau>0\}$.) Since $\tau=-\infty$ on $X_{sing}$, the set  
\[
		\Omega=\{x\in X\colon \rho(x)<M_3,\ \tau(x)> 0\}
\]
is contained in the regular part of $X$. By a small perturbation 
one can in addition achieve that $0$ is a regular value of $\tau$, 
$M$ is a regular value of $\rho$,  and the level sets 
$\{\rho=M\}$ and $\{\tau=0\}$ intersect transversely. 
Denote their intersection manifold by $\Sigma$.
There is a neighborhood $U_\Sigma$ of $\Sigma$ in $X$
with $\overline U_\Sigma \subset \{\rho>M_2\} \cap X_{reg}$.

We are now in the same geometric situation as in \cite[Subsection 6.5]{ACTA}.
(See especially the proof of Lemma 6.9 in \cite{ACTA}. The fact that 
our $X$ is not necessarily a manifold is unimportant since
$\overline \Omega \subset X_{reg}$.)  For $s\in [0,1]$ set  
\[
		\rho_s =(1-s)\tau + s(\rho-M),\quad G_s=\{\rho_s<0\} \cap \{\rho< M_3\}.
\]
The Levi form of $\rho_s$, being a convex combination of the
Levi forms of $\tau$ and $\rho$, is positive on the linear span of 
the eigenspaces corresponding to the positive eigenvalues of the 
Levi form of $\rho$. 
Therefore $G_s$ is strongly $(n-1)$-convex at each smooth
boundary point for every $s\in [0,1]$. As the parameter $s$ increases from 
$s=0$ to $s=1$, the domains $G_s\cap \{\rho<M\}$ increase from 
$\{\tau < 0,\ \rho<M\}$ to $G_1=\{\rho < M\}$. 
(The sets $G_s\cap \{M < \rho <M_3\}$ decrease with $s$, 
but that part will not be used.) All hypersurfaces $\{\rho_s = 0\} = bG_s$ 
intersect along $\Sigma$.  Since $d\rho_s = (1 - s)d\tau + s d\rho$ and
the differentials $d\tau$ and $d\rho$ are linearly independent along $\Sigma$,
each hypersurface $bG_s$ is smooth near $\Sigma$.
By a generic choice of $\rho$ and $\tau$ we can insure that only for finitely many values 
of $s \in [0,1]$ does the critical point equation $d\rho_s=0$ have a solution 
on $bG_s\cap \Omega$, and in this case there is exactly one solution. 
Therefore $bG_s$ has non smooth points only for finitely many values
of $s\in [0,1]$.

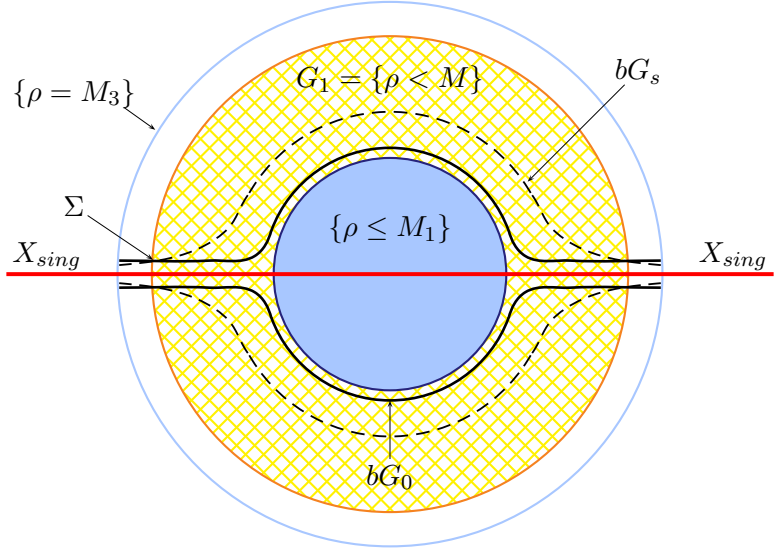
\begin{figure}[ht]
\psset{unit=0.6 cm} 
\begin{pspicture}(-9,-6.1)(9,6.1)

%
%

\pscircle[linecolor=myblue](0,0){6.06}       
\pscircle[linecolor=OrangeRed,fillstyle=crosshatch,hatchcolor=yellow](0,0){5.3}       

%
%
\pscircle[linecolor=DarkBlue,fillstyle=solid,fillcolor=myblue](0,0){2.6}     
    	            
\psline[linecolor=red,linewidth=1.5pt](-8.5,0)(8.5,0)                                  
\rput(-7.6,0.4){$X_{sing}$}
\rput(7.6,0.4){$X_{sing}$}

%
%
\psarc[linewidth=1pt](0,0){2.8}{20}{160}
\psarc[linewidth=1pt](0,0){2.8}{200}{340}

%
%
%
\psecurve[linewidth=1pt](3,5)(2.63,0.97)(3,0.4)(4,0.3)(6,0.3)(7,0.3) 
\psecurve[linewidth=1pt](3,-5)(2.63,-0.97)(3,-0.4)(4,-0.3)(6,-0.3)(7,-0.3) 
\psecurve[linewidth=1pt](-3,5)(-2.63,0.97)(-3,0.4)(-4,0.3)(-6,0.3)(-7,0.3) 
\psecurve[linewidth=1pt](-3,-5)(-2.63,-0.97)(-3,-0.4)(-4,-0.3)(-6,-0.3)(-7,-0.3)

%
%
%
\psarc[linestyle=dashed,linewidth=0.7pt](0,0){3.6}{30}{150}
\psarc[linestyle=dashed,linewidth=0.7pt](0,0){3.6}{210}{330}

%
%
%
\psecurve[linestyle=dashed,linewidth=0.7pt](2,3)(3.15,1.75)(3.6,1)(4.4,0.5)(6,0.2)(7,0.2)
\psecurve[linestyle=dashed,linewidth=0.7pt](2,-3)(3.15,-1.75)(3.6,-1)(4.4,-0.5)(6,-0.2)(7,-0.2)
\psecurve[linestyle=dashed,linewidth=0.7pt](-2,3)(-3.15,1.75)(-3.6,1)(-4.4,0.5)(-6,0.2)(-7,0.2)
\psecurve[linestyle=dashed,linewidth=0.7pt](-2,-3)(-3.15,-1.75)(-3.6,-1)(-4.4,-0.5)(-6,-0.2)(-7,-0.2)

\rput(0,4.3){$G_1=\{\rho<M\}$}						                            

\psline[linewidth=0.2pt]{<-}(3.05,2.05)(5.2,4.2)
\rput(5.5,4.5){$bG_s$}

\rput(0,1){$\{\rho\le M_1\}$}	

\rput(-7,1.5){$\Sigma$}
\psline[linewidth=0.2pt]{->}(-6.7,1.3)(-5.25,0.32)

\rput(-7,4){$\{\rho=M_3\}$}
\psline[linewidth=0.2pt]{->}(-5.9,3.7)(-5.2,3.2)

\rput(0,-4.5){$bG_0$}  
\psline[linewidth=0.2pt]{->}(0,-4.1)(0,-2.8)

\end{pspicture}
\caption{The sets $G_s$.} 
\label{Fig4}
\end{figure}

Fix two values of the parameter, say $0 \le s_0 < s_1 \le 1$. 
Consider first the {\em noncritical case} when $d\rho_s \ne 0$ on
$bG_s\cap \Omega$ for all $s\in[s_0,s_1]$, and hence all boundaries 
$bG_s$ for $s\in[s_0, s_1]$ are smooth.
By attaching to $G_{s_0}$ finitely many small 2-convex bumps of the type 
used in Lemma \ref{smallstep} and contained in $G_1\cup U_\Sigma$ we cover 
the set $G_{s_1}\cap \Omega$. (See \cite[p.\ 180]{ACTA} for a more detailed
description.) Using Lemma \ref{smallstep} at each bump 
we push the boundary of the central map in the spray outside the bump 
while keeping control on the compact subset $\overline U \subset D$. 
After a finite number of steps the boundary of the central map 
lies outside $G_{s_1}\cap \Omega$ and inside $G_1\cup U_\Sigma$.
In the sequel this will be called the {\em noncritical procedure}.

It remains to consider the values $s\in [0, 1]$ for which $bG_{s}$ has a
non smooth point (the {\em critical case}). 
We begin by discussing the most difficult case $\dim X=2$
when there is least space to avoid the critical points.
The functions $\rho$ and $\tau$ are then $1$-convex and hence \spsh.
As in \cite[p.\ 180]{ACTA} we introduce the function 
\[
	h(x)=\frac{\tau(x)}{\tau(x)+ M - \rho(x)}, \quad x\in\Omega.
\]
%
%
A generic choice of $\tau$ insures that $h$ is a Morse function. 
Note that $\{h=s\} = \{\rho_s=0\} =bG_s$.
The critical points of $h$ coincide with critical points of $\rho_s$ on 
$\{\rho_s=0\}$, and the Levi form of $h$ at a critical point is positive definite
\cite[p.\ 180]{ACTA}.

To push the boundary over a critical level of $h$ we shall apply 
Lemma 6.7 in \cite[p.\ 177]{ACTA} (see also \cite[\S 4]{FK}).
Let $p$ be a critical point of $h$, with $h(p)=c \in (0,1)$.
(Our $h$ corresponds to $\rho$ in \cite{ACTA}.)
It suffices to consider the case when the Morse index of $p$ is 
either 1 or 2 since we cannot approach a minimum of $h$ by the
noncritical procedure. Choose a neighborhood $W \subset X$ of $p$ on 
which $h$ is \spsh. Lemma 6.7 in \cite{ACTA} furnishes a new function $\wt h$ 
(denoted $\tau$ in \cite{ACTA})
which is \spsh\ on $W$, while outside of $W$ each level set 
$\{\wt h=\e\}$ (for values $\e$ close to $0$) coincides with 
a certain level set $\{h=c(\e)\}$, such that $\wt h$
satisfies the following properties (see fig.\ \ref{Fig5}). 
The sublevel set $\{\wt h\le 0\}$ is contained in the union of the sublevel set 
$\{h \le c_0\}$ for some $c_0<c$ (close to $c$) and a totally real disc $E$
(the unstable manifold of the critical point $p$ with respect 
to the gradient flow of $h$). Furthermore, for a small $d>0$ 
with $c_0 < c-d$ we have 
\begin{equation}
\label{inclusions}
		\{h \le c+d\} \subset \{\wt h \le 2d\} \subset \{h< c+3d\},
\end{equation}
$\wt h$ has no critical values on $(0,3d)$,
and $h$ has no critical values on $[c-d,c+3d]$ except for $h(p)=c$.

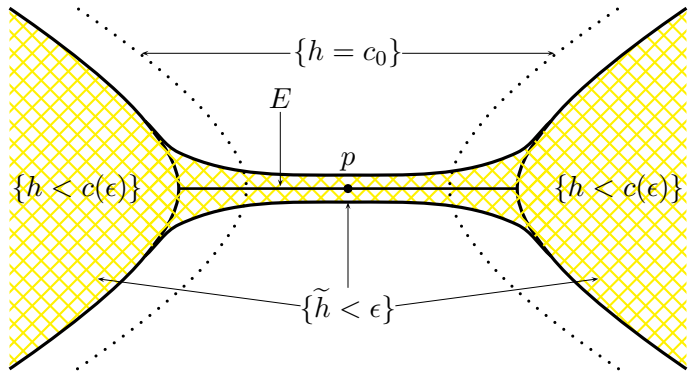
\begin{figure}[ht]
\psset{unit=0.6cm, xunit=1.5, linewidth=0.7pt} 
\begin{pspicture}(-4.5,-4.3)(4.5,4.3)

%
%
%
%

\pscustom[fillstyle=crosshatch,hatchcolor=yellow,linewidth=1.2pt]
{\pscurve[liftpen=1](5,4)(3,1.5)(2.5,0)(3,-1.5)(5,-4)                 
}

\pscustom[fillstyle=crosshatch,hatchcolor=yellow,linewidth=1.2pt]
{
\pscurve[liftpen=1](-5,4)(-3,1.5)(-2.5,0)(-3,-1.5)(-5,-4)             
}

\pscustom[fillstyle=crosshatch,hatchcolor=yellow,linestyle=none]        
{
\pscurve(-3,-1.5)(-2.5,-0.8)(0,-0.3)(2.5,-0.8)(3,-1.5) 
\psline[linestyle=dashed,linewidth=0.2pt](3,-1.5)(3,1.5)
\pscurve[liftpen=1](3,1.5)(2.5,0.8)(0,0.3)(-2.5,0.8)(-3,1.5)
\psline[linestyle=dashed,linewidth=0.2pt](-3,1.5)(-3,-1.5)
}

%
%
%
\psline[linewidth=1pt](-2.5,0)(2.5,0)                                                  
\psecurve[linewidth=1.2pt](5,4)(3,1.5)(2.5,0.8)(0,0.3)(-2.5,0.8)(-3,1.5)(-5,4)         
\psecurve[linewidth=1.2pt](5,-4)(3,-1.5)(2.5,-0.8)(0,-0.3)(-2.5,-0.8)(-3,-1.5)(-5,-4)  

%
%
%
\pscurve[linestyle=dotted,linewidth=1.2pt](4,4)(2,1.5)(1.5,0)(2,-1.5)(4,-4)                  
\pscurve[linestyle=dotted,linewidth=1.2pt](-4,4)(-2,1.5)(-1.5,0)(-2,-1.5)(-4,-4)            

\rput(0,3){$\{h=c_0\}$}
\psline[linewidth=0.2pt]{->}(0.8,3)(3.05,3)
\psline[linewidth=0.2pt]{->}(-0.85,3)(-3.05,3)

%
%
%
%
\rput(4,0){$\{h< c(\epsilon)\}$}                                           
\rput(-4,0){$\{h<c(\epsilon)\}$}

\psline[linewidth=0.2pt]{->}(0,-2.2)(0,-0.3)
\psline[linewidth=0.2pt]{<-}(-3.7,-2)(-0.7,-2.6)
\psline[linewidth=0.2pt]{->}(0.7,-2.6)(3.7,-2)
\rput(0,-2.6){$\{\wt h < \epsilon\}$}

\psline[linewidth=0.2pt]{->}(-1,1.7)(-1,0.05)
\rput(-1,2){$E$}

\psdot(0,0)
\rput(0,0.6){$p$}

\end{pspicture}
\caption{The level sets of $\widetilde h$}
\label{Fig5}
\end{figure}

By the noncritical procedure applied with the function 
$h$ we push the boundary 
of the central map of the spray into the set $\{c-d < h < c\}$. 
Let $\wt f$ denote the new spray. For parameters $t\in\C^N$ 
sufficiently close to $0$ the map $\wt f_t$ also has boundary values in 
$\{c-d <h <c\}$. Since $\dim_\R E\le 2$, we can find $t$ arbitrarily close 
to $0$ such that $\wt f_t(bD) \cap E= \emptyset$. 
By translation in the $t$ variable we can choose 
$\wt f_t$ as the new central map of the spray. 

Since $\{\wt h\le 0\}\subset \{h \le c_0\}\cup E \subset  \{h \le c-d\}\cup E$,
the above insures that $\wt h >0$ on $\wt f_t(bD)$.
Since $\wt h$ has no critical values on $(0,3d)$,
we can use the noncritical procedure with $\wt h$ to push the 
boundary of the central map into the set $\{\wt h > 2d\}$, appealing to
Lemma \ref{smallstep}. As $\{\wt h > 2d\} \subset \{h > c+d\}$ by (\ref{inclusions}),
we have thus pushed the image of $bD$ across the critical level 
$\{h=c\}$ and avoided running into the critical point $p$. 
Now we continue with the noncritical procedure applied 
with $h$ to reach the next critical level of $h$. 

This concludes the proof for $n=2$.  The same procedure can be adapted
to the case $n=\dim_\C X >2$ by considering the appropriate two dimensional 
slices on which the function $\rho$ is strongly plurisubharmonic. 
Alternatively, we can apply the same geometric construction 
as in \cite{BD3} to keep the boundary of the central map at a positive 
distance from the critical points of $\rho$.
\end{proof}

%
%
%
%
\textit{Proof of Theorem \ref{main}.}
Let $d$ denote a complete distance function on $X$.
We denote the initial map in Theorem \ref{main} 
by $f_0\colon \bar D\to X$. By Theorem \ref{approximation2}
we may assume that $f_0$ is holomorphic in a neighborhood of $\bar D$
in an open Riemann surface $S\supset \bar D$ and 
$f_0(bD)\subset (X_c)_{reg}$. Here $X_c=\{\rho>c\}$ is the set
on which $\rho$ is assumed to have at least two positive eigenvalues.

Choose an open relatively compact subset $U\ss D$ and a number $\e>0$.
It suffices to find a proper holomorphic map 
$g\colon D\to X$ such that $\sup_{z\in U} d(f_0(z),g(z)) < \e$
and such that $g$ agrees with $f_0$ to order $k$ at each of the 
given points $z_j \in D$; a sequence of proper maps 
$g_\nu$ as in Theorem \ref{main} is then obtained by 
Cantor's diagonal process.

Let $\sigma$ denote the union of $\{z\in D\colon f_0(z) \in X_{sing}\}$
and the finite set $\{z_j\}\subset D$ on which we wish to interpolate
to order $k\in \N$; thus $\sigma$ is a finite subset of $D$.
Lemma \ref{sprays-exist} furnishes a spray of maps 
$f\colon \bar D\times P\to X$ of class $\cA^2(D)$, with the
given central map $f_0$ and the exceptional set $\sigma$ of order $k$, 
such that $f_t(bD) \subset (X_c)_{reg}$ for each $t\in P \subset\C^N$.

Set $f^0=f$, $c=c_0$, and choose an open
subset $P_0\ss P$ containing the origin $0\in \C^N$. 
Choose a number $c_1>c_0$ such that 
$c_0< \rho(f^0_t(z)) <c_1$ for all $z\in bD$ and $t\in P_0$,
and then choose an open subset $U_0\Subset D$ containing 
$\sigma\cup U$ such that 
$f^0_t(\bar D\bs U_0)\subset \{x\in X\colon c_0<\rho(x)<c_1\}$
for all $t\in P_0$. Choose a sequence $c_0 < c_1< c_2\cdots$ 
with the given initial numbers $c_0$ and $c_1$ such that 
$\lim_{j\to\infty} c_j=+\infty$. 
Also choose a decreasing sequence $\e_j>0$ with $0<\e_1<\e$ 
such that for each $j\in \N$ we have
\[
	\bigl( x,y\in X,\ \rho(x)<c_{j+1},\ d(x,y)<\e_j\bigr)
	 \Rightarrow |\rho(x)-\rho(y)|<1.
\]	 

We shall inductively find a sequence of sprays 
$f^j\colon \bar D\times P_j\to X$ of class $\cA^2(D)$
with the exceptional set $\sigma$ of order $k$, 
with $P=P_0\supset P_1\supset P_2\supset\cdots$, 
and a sequence of open sets 
$U_0\subset U_1\subset \cdots\subset \cup_{j=1}^\infty U_j = D$
satisfying the following properties for each $j\in\Z_+$
and $t\in P_j$:
\begin{itemize}
\item[(i)]   $f^j_t(bD) \subset \{x\in X_{reg} \colon c_j< \rho(x) < c_{j+1}\}$,
\item[(ii)]  $f^j_t(\bar D\bs U_j) \subset \{x\in X\colon c_j< \rho(x) < c_{j+1}\}$,
\item[(iii)] $f^j_t(\bar D\bs U_{j-1}) \subset \{x\in X\colon c_{j-1}< \rho(x) < c_{j+1} \}$,
\item[(iv)] $d\bigl(f^{j}_0(z),f^{j-1}_0(z)\bigr) <\e_j 2^{-j}$ for $z\in U_{j-1}$, and
\item[(v)]  $f^{j}_0$ and $f^{j-1}_0$ are homotopic, and they have the same 
$k$-jets at each of the points in $\sigma$.
\end{itemize}

For $j=0$ the properties (i) and (ii) hold while the remaining properties
are vacuous. (In (iii) we take $U_{-1}=U_0$ and $c_{-1}=c_0$.)
Assuming that we already have sprays $f^0,\ldots, f^j$ satisfying these 
properties, Proposition \ref{bigstep} applied to $f=f^j$ furnishes a new spray 
$f^{j+1}$ (called $g$ in the statement of that Proposition) 
satisfying (i), (iii), (iv) and (v). Choose an open set $U_{j+1} \Subset D$  
with $U_j\subset U_{j+1}$ such that (ii) holds (this is possible by continuity
since (i) already holds and we are allowed to shrink the parameter set $P_{j+1}$). 
Hence the induction proceeds. When choosing the sets $U_j$ we can easily 
insure that they exhaust $D$.

Conditions (i)--(v) imply that the sequence of central maps 
$f^j_0\colon \bar D\to X$ $(j\in\Z_+)$ converges uniformly on compacts in $D$ 
to a proper holomorphic map $g\colon D\to X$ satisfying 
$d(f_0(z),g(z)) < \e$ $(z\in \overline U_0)$ and such that the 
$k$-jet of $g$ agrees with the $k$-jet of $f_0$ at every point of $\sigma$.
In addition,  we can combine the homotopies 
from $f^j_0$ to $f^{j+1}_0$ $(j=0,1,\ldots)$ 
to obtain a homotopy from $f_0|_D$ to $g$. 
This completes the proof of Theorem \ref{main}.

\section{Appendix: Approximation of holomorphic vector subbundles}
In the proof of  Lemma \ref{transition-map} we used the following 
approximation result:  

\begin{theorem} 
{\rm (Heunemann \cite{Heunemann1})}
If $D$ is a relatively compact strongly pseudoconvex domain in a Stein 
manifold $S$ and $E\subset \bar D\times \C^n$ is a continuous 
complex vector subbundle of the trivial bundle over $\bar D$
such that $E$ is holomorphic over $D$ then $E$ can be uniformly approximated 
by holomorphic vector subbundles $\wt E \subset U\times \C^n$ 
over small open neighborhoods $U\subset S$ of $\bar D$.
\end{theorem}

We offer a simple proof of this useful result. Choose a complementary 
to $E$ subbundle $G\subset \bar D\times \C^n$ of the same class 
$\cA(D)$ (the existence of such $G$ follows from Cartan's Theorem B
for vector bundles of class $\cA(D)$ \cite{Heunemann2}, \cite{Leiterer2}).
Let $\Pi \colon \bar D \times\C^n \to E$  denote the 
fiberwise $\C$-linear projection with kernel $G$ and image $E$.
By the Oka-Weil theorem we approximate $\Pi$ uniformly on $\bar D$ by a 
holomorphic fiberwise linear map $\Pi' \colon U'\times \C^n \to U'\times \C^n$
over an open set $U' \supset \bar D$. In general $\Pi'$ will fail to be 
a projection map on the fibers, but this can be corrected by the following 
simple device (see e.g.\ \cite{GLR}):

\smallskip
{\em Let $C$ be a positively oriented simple closed curve in $\C$,
and let $L\in {\rm Lin}_\C(\C^n,\C^n)$ be a linear map with no
eigenvalues on $C$. Then $\C^n= V_+\oplus V_-$, where $V_+$ resp.\
$V_-$ are $L$-invariant subspaces of $\C^n$ spanned by the generalized
eigenvectors of $L$ corresponding to the eigenvalues 
inside resp.\ outside of $C$. The map
\[
	\cP(L)= \frac{1}{2\pi i} \int_C \left(\zeta I-L\right)^{-1}\, d\zeta
\]
is a projection onto $V_+$ with kernel $V_-$.}
\smallskip

Choose a curve $C \subset \C$ which encircles $1$ but not $0$;
for instance, $C=\{\zeta \in \C \colon |\zeta -1|=1/2\}$. 
Let $\cP$ denote the associated projection operator. 
If $L\in {\rm Lin}_\C(\C^n,\C^n)$ is a projection then $\cP(L)=L$. 
If $L'$ is  near a projection $L$ then each eigenvalue 
of $L'$ is either near $0$ or near $1$, and hence $\cP(L')$ is a 
projection which is close to $L$ and has the same rank as $L$.

Assuming that $\Pi'$ is sufficiently close to $\Pi$ on $\bar D$
it follows that for each point $z$ in an open set $U'$ 
with $\bar D\subset U\subset U'$
the map $\wt \Pi_z = \cP(\Pi'_z) \in {\rm Lin}_\C(\C^n,\C^n)$ is 
a projection of the same rank as $\Pi_z$ and it depends 
holomorphically on $z\in U$. The map $\wt \Pi\colon U\times \C^n \to  U\times\C^n$ 
with fibers $\wt\Pi_z$ is then a projection onto a holomorphic vector 
subbundle $\wt E \subset  U\times\C^n$ whose restriction to  $\bar D$ 
is uniformly close to $E$, and $\wt G = {\rm ker}\, \wt \Pi$
is a holomorphic vector subbundle of $U\times\C^n$ 
whose restriction to $\bar D$ is uniformly close to $G$.

%
%
%
%
%

\medskip
\textit{Acknowledgements.} 
The first named author wishes to thank the Laboratoire de Math\'ematiques
E.\ Picard, Universit\'e Paul Sabatier de Toulouse, for its hospitality, 
and the EGIDE program for support during a part of the work.
The second named author thanks D.\ Barlet, M.\ Brunella, J.-P.\ Demailly, 
C.\ Laurent-Thi\'ebaut, J.\ Leiterer, I.\ Lieb, J.\ Michel, M.\ Range, 
N.\ \O vrelid and J.-P.\ Rosay for helpful discussions, 
and the Institut Fourier, Universit\'e de Grenoble, for support and hospitality 
during a part of the work. We also thank the referee for pertinent remarks.

\bibliographystyle{amsplain}

\end{document}